\documentclass[11pt]{article}
\usepackage{amsmath,latexsym,amssymb,amsfonts,amsbsy, amsthm}
\usepackage [latin1]{inputenc}
\usepackage{color, xcolor}

\usepackage{graphicx}
\usepackage{soul, color}
\usepackage[margin = 1.5in]{geometry}
\newtheorem{theorem}{Theorem}[section]

\newtheorem{corollary}[theorem]{Corollary}

\newtheorem{lemma}[theorem]{Lemma}

\newtheorem{remark}[theorem]{Remark}

\def\mau{\mathbf{u}}

\def\dv{\mathrm{div}}
\let\f=\frac

\let\va=\varepsilon

\let\f=\frac

\def\dv{\mathrm{div}}

\def\curl{\mathop{\rm curl}\nolimits}

\begin{document}

\title{\textbf{Zero Viscosity Limit of Steady Compressible Shear Flow with Navier-Slip Boundary }}

\author{ \bf{ Wenbin Li }\footnote{liwenbin@seu.edu.cn.\ School of Mathematics, Southeast University, Nanjing, China.}\hspace{10 mm} \bf{Chunhui Zhou}\footnote{zhouchunhui@seu.edu.cn. School of Mathematics, Southeast University, Nanjing, China.} }

\date{}
\maketitle

\begin{abstract}
We investigate the existence and the zero viscosity limit of   steady compressible shear flow with Navier-slip boundary condition in the absence of any external force in a two-dimension domain $\Omega=(0,L)\times(0,2)$. More precisely, under the assumption that the Mach number $\eta<\va^{\f12+}$ and $L\ll1$, we prove the existence of smooth solutions  to steady compressible Naiver-Stokes equations near plane Poiseuille-Couette flow as well as the convergence of the solutions obtained above to the solutions of steady incompressible Euler equations when the viscous $\va$ tends to zero.
\end{abstract}

\bigskip
\renewcommand{\theequation}{\thesection.\arabic{equation}}
\setcounter{equation}{0}
%%%%%%%%%%%%%%%%%%%%%%%%%%%%%%%%%%%%%%%%%%%%%%

\section{Introduction}

In this paper, we shall study the existence and the zero viscosity limit of   steady compressible shear flow near plane Poiseuille-Couette flow in the absence of any external force in a 2-dimension domain $\Omega=(0,L)\times (0,2)$. Considering the dimensionless spatial coordinates $x/L^*,\ y/L^*$, the dimensionless velocity and density $\mau/V^*,\ \rho/\rho^*$, here $\mau=(u,v)$, steady compressible Navier-Stokes equations can be written into the following dimensionless form:
\begin{eqnarray}
\dv(\rho\mau)=0,\label{0.3}\\
\rho\mau\cdot\nabla u-\mu\va\Delta u-\lambda\va\partial_x\dv\mau+\f1{\eta^2}\partial_xP=0,\label{0.4}\\
\rho\mau\cdot\nabla v-\mu\va\Delta v-\lambda\va\partial_y\dv\mau+\f1{\eta^2}\partial_yP=0,\label{0.5}
\end{eqnarray}
here $\va=1/Re,\ Re=\f{\rho^*V^*L^*}{\mu^*}$ is the Reynolds number; $P$ is the pressure for isentropic flows  given by
$P(\rho)=a\rho^\gamma$ with $a$ being a positive constant and
$\gamma>1$ being the specific heat ratio,  $\mu>0,\mu'>0$ are the scaled shear and bulk viscosity, $\lambda=\mu'+\f13\mu$; $\eta=\frac{V^*}{c}$ is the Mach number and $c=\sqrt{p'(\rho^*)}$  the speed of sound. To simplify the presentation, we will take $a=\mu=\lambda=1$ in the following of this paper.

Before introducing the main results, let us briefly review the related known results. The vanishing viscous limit
for unsteady compressible Naiver-Stokes equations with slip or non-slip boundary conditions in Sobolev space or in analytic setting has been extensively studied  in cf.\cite{LW,P,W,WWZ,WW,WXY,XY} and the reference therein.  Recently, in \cite{ADM,ZZ}, the authors showed the linear  stability  of the Couette flow in compressible fluid in 2D and 3D. More recently, the authors in \cite{MWWZ,YZ1,YZ}  showed the linear spectral instability analysis of planar Poiseuille-type flow and strong boundary layer in compressible fluid.

In the steady setting, the existence of weak solutions to steady compressible Navier-Stokes equations with homogeneous  boundary condition for the velocity and arbitrary  large external force have been proved in \cite{L-2,J-Zhou,B-N,P-W} and the reference therein. The strong solutions of the steady  compressible Navier-Stokes equations  with
homogenous boundary conditions  have been
investigated in cf. \cite{C-J,D-J-J-Y,J-K,N-N-P,Pa,V,V-Z} and the reference therein. To the best of our knowledge, the question about the correct formulation and the existence of weak solutions with inhomogeneous boundary conditions is still an open problem. Concerning the strong solutions
with inhomogeneous boundary data,  the authors in \cite{K-K,k-k-2,K-K-3} proved the existence of strong solutions to the
stationary problems with an inflow boundary condition under the assumption that the Reynolds number of the flow is small.
 In \cite{P-R-S,P-R-S-1} the authors study the shape optimization problem for
 the compressible Navier-Stokes equations with inhomogeneous boundary value in the case of small
Reynolds and Mach numbers. In \cite{PT,P-P}, the authors investigate the existence of strong
solutions with inhomogeneous slip boundary conditions on the velocity
under the assumption that the friction coefficient on the boundary is large enough. In \cite{GJZ}, the authors proved the existence of strong solutions to steady compressible Navier-Stokes equations with inflow boundary condition around the constant flow $\tilde{\mathbf{U}}=(1,0)$ in a channel without the smallness assumption of the Reynolds number and the Mach number.

In this paper, we will study the structure stability of steady compressible shear flow near plane Poiseuille-Couette flow $\mathbf{U}^0=(u_s,0)$  in $\Omega$  in the absence of any external force.
For $\alpha_0>0,\alpha_1\geq0,\alpha_2\geq0,\alpha_1+\alpha_2>0$,  $u_s$ is defined as following:
$$u_s=\alpha_0+\alpha_1y+\alpha_2y(2-y),\ 0<y<2.$$
If we take $P^0\equiv1$, then it is easy to check that the triple $(u_s,0,1)$
satisfies the stationary incompressible Euler equations:
\begin{equation} \label{main.euler}
\begin{cases}
\bold{U}^0 \cdot \nabla \bold{U}^0 + \nabla P^0 = 0 \\
\nabla \cdot \bold{U}^0 = 0 \\
\bold{U}^0 \cdot \bold{n}|_{y = 0, y = 2} = 0.
\end{cases}
\end{equation}
Here  $\bold{n}$ denotes the unit normal vector on $\{y = 0\}$ and $\{y = 2\}$.

On the other hand, if we take $P^{NS}=1-2\va\alpha_2x$, then  the triple $(u_s,0,P^{NS})$
satisfies the stationary incompressible Navier-Stokes equations:
\begin{equation} \label{main.NS}
\begin{cases}
\bold{U}^0 \cdot \nabla \bold{U}^0-\va\Delta \bold{U}^0 + \nabla P^{NS} = 0 \\
\nabla \cdot \bold{U}^0 = 0.
\end{cases}
\end{equation}
We are interested in the compressible perturbation of the stationary incompressible Navier-Stokes equations (\ref{main.NS}) as well as the zero viscosity limit of steady slightly compressible flow to steady incompressible Euler flow (\ref{main.euler}).
The Navier-Stokes equations for a steady isentropic compressible viscous
flow is a mixed system of hyperbolic-elliptic type, as the momentum
equations are an elliptic system in the velocity, while the continuity
equation is hyperbolic in the density. Therefore,  it is necessary to prescribe the density on
the part of inflow boundary ($\mathbf{U}^0\cdot\mathbf{n}<0$), where $\mathbf{n}$ is the
outward unit normal to the boundary.  For the Navier-Stokes equations in a channel,  there are no natural boundary conditions at $\{x=0\},\{x=L\}$,  thus part of the mathematical challenge is to impose proper boundary conditions for the unknowns, especially the boundary condition for the density on the inflow part of the bounary, to ensure its solvability as well as the uniform-in-$\va$ estimates. In our method, in the process of the zero viscosity limit $\va\rightarrow0$ of steady compressible flow, we need to prescribe the value of the pressure gradient in the flow direction, i.e. $P_x$, on the inflow boundary to control  the convection terms. Fortunately, we find from the second equation in the momentum equations that pressure gradient in the vertical direction, i.e. $P_y$, on the inflow boundary can also be well controlled if we can control the curl of the flow on the inflow boundary. Our boundary condition is described as following:
\begin{eqnarray}
&&\mau^\va\cdot\mathbf{n}|_{y=0,2}=0,\ u_y^\va(x,0)=u_{sy}(0)+b_{0}(x),\ u_y^\va(x,2)=u_{sy}(2)+b_{1}(x),\label{0.b1}\\
&&u^\va|_{x=0}=u_s+a_{1}(y),\ v_x^\va|_{x=0}=a_{2}(y),\\
&&v^\va|_{x=L}=a_{3}(y),\ \partial_x(u_y^\va-v^\va_x)|_{x=L}=a_{4}(y),\\
&&\rho_x^\va|_{x=0}=-\f2\gamma\va\eta^2\alpha_2+h_0(y),\ \rho^\va(0,0)=1.\label{0.b2}
\end{eqnarray}
with compatibility condition:
\begin{eqnarray}
&&a_1'(0)=b_0(0),\ a_1'(2)=b_1(0),\ a_4(0)=b_0'(L),\ a_4(2)=b_1'(L),\nonumber\\
&&a_2(0)=a_2(2)=a_3(0)=a_3(2)=a_4''(0)=a_4''(2)=0.\label{0.2}
\end{eqnarray}
It is easy to check that the under the slightly compressible perturbation of the stationary incompressible flow with $P=\rho^\gamma$, although the strict parrel flow $(u,v,\rho)=(u_s(y),0,(1-2\va\eta^2\alpha_2x)^{\f1\gamma})$ satisfy the momentum equations (\ref{0.4}) and (\ref{0.5}), the mass equation (\ref{0.3}) fails as we will have $\dv(\rho\mau)=-\f2\gamma\alpha_2\eta^2\va u_s(1-2\va\eta^2\alpha_2x)^{\f1\gamma-1}$. In fact, even if $b_0(x)=b_1(x)=a_i(y)=h_0(y)\equiv0,$  $ i=1,2,3,4$, the solution $(\rho^\va,\mau^\va)$ to system (\ref{0.3})-(\ref{0.5}) with boundary condition (\ref{0.b1})-(\ref{0.b2}) is nontrivial.

Now let us introduce the notations used throughout this paper.
\\[1mm]
{\sc Notation:} \ Let $G$ be an open domain in $\mathbb{R}^N$.
We denote by $L^p(G)$ ($p\geq 1$) the Lebesgue spaces, by $W^{s,p}(G)$ ($p\geq 1$) the Sobolev spaces with $s$ being a real number,
by $H^k(G)$ ($k\in \mathbb{N}$) the Sobolev spaces $W^{k,p}(G)$ with $p=2$, and
by $C^k(G)$ (resp. $C^k(\overline{G})$) the space of $k$th-times continuously differentiable functions in $G$ (resp. $\overline{G}$).
 We use $|\cdot|_{k,p}$ to denote the standard norm in $W^{k,p}(G)$ at the boundary $\partial\Omega$ and $|\cdot|$ for the norm in $L^2(G)$   throughout this paper.
 $\|\cdot\|_{k,p}$ stands for the standard norm in $W^{k,p}(G)$ and $\|\cdot\|$ for the norm in $L^2(G)$.
 We also use $|\cdot|_{L^\infty}$ to denote $|\cdot |_{L^\infty(\Omega)}=\text{ess~sup}_{\Omega}|\cdot |$.
 The symbol $\lesssim$ means that the left side is less than the right side multiplied by some constant.
Let $a\in {\mathbb R}$ be a real number, then $a+$ means any real number a little bigger than $a$, $a-$ means any real number a little smaller than $a$.
 For $G=\Omega$, we define
\begin{eqnarray}
|u|_{L^\infty_xL^2_y}=\text{ess}\sup_{0<x<L}|u(x,\cdot)|_{L^2(0,2)}.\nonumber\end{eqnarray}
 %%%%%%%%%%%%%%%%%%%%%%

We also define a smooth cut-off function $\chi(t)\in C^\infty([0,\infty))$ satisfying $|\chi|\leq1,\ |\chi|_{C^4}\leq C$ and
\begin{equation}\label{cutoff}\chi(t)=\begin{cases}1,\ 0\leq t\leq 1/2,\\0,\ t\geq1,\end{cases}\end{equation}
here $C>0$ is a finite constant.

Next we define the space $\mathcal{A}$  by
\begin{eqnarray}
\mathcal{A}=&& \{\rho\in H^2(\Omega)\ |\ \rho_x|_{x=0}=0,\ \rho(0,0)=0\ \|\rho\|_{\mathcal{A}}<\infty\},\label{a}
\end{eqnarray}
with the norm
\begin{eqnarray}\|\rho\|_{\mathcal{A}}=&&\|\rho\|_{H^1}+\va^{\f12}\eta|\nabla\rho|_{L^\infty_xL^2_y}+\va\eta\|\rho_{xx}\|+\va^{\f12}\|\nabla\rho_{y}\|
+\va^{\f32}\eta^2|\rho_{xx}|_{L^\infty_xL^2_y}\nonumber\\
&&+\va\eta|\nabla\rho_{y}|_{L^\infty_xL^2_y}+|\rho(0,\cdot)|_{H^1}+\va^{\f12}|\rho_{yy}(0,\cdot)|\end{eqnarray}
and for $2\leq p<\infty$, we define the space $\mathcal{B}$:
\begin{eqnarray}
\mathcal{B}=&& \{\mau=(u,v)\in W^{2,p}(\Omega)\times  W^{2,p}(\Omega), \;\curl\mau\in W^{2,p}(\Omega),\ (L-x)\mau\in H^3(\Omega)\nonumber\\
&&|\ \|\mau\|_{\mathcal{B}}<\infty,\ u_y|_{y=0,2}=v|_{y=0,2}=u|_{x=0}=v_{x}|_{x=0}=v|_{x=L}=\curl\mau_x|_{x=L}=0\},\nonumber\\\label{b}
\end{eqnarray}
where the norm $\|\cdot\|_{\mathcal{B}}$ is defined by
\begin{eqnarray*}\|\mathbf{u}\|_{\mathcal{B}}:=&&\|\mau\|_{H^1}+\va^{\f12}\|\nabla^2\mau\|+\va^{\f12}\|u_y-v_x\|_{H^1}+\va^{2-\f2p}\|\mau\|_{2,p}+\va^{2-\f2p}\|u_y-v_x\|_{W^{2,p}}\\
&&+\va^{\f32}|\dv\mau_{yy}(0,\cdot)|+\va^2\eta\|\dv\mau_{xx}\|+\va^{\f32}(\|(L-x)\nabla^2u_y\|\nonumber\\
&&+\|(L-x)\nabla^3v\|)+\va^2\eta\|(L-x)u_{xxx}\|+\va^{\f32}\|\nabla\dv\mau_y\| \end{eqnarray*}
Finally, we denote by
\begin{equation}
\Lambda=\eta^{-2}|h_0|_{H^3}+\sum_{i=1}^4|a_i|_{H^4}+|b_0|_{H^4}+|b_1|_{H^4}\label{0.1}
\end{equation}

Our  main result reads as follows.
%%%%%%%%%%%%%%%%%%%%%%%%%%%%%%%%%%%%%%%%%%%%%%%
\begin{theorem}\label{th.main}
Let $0<L\ll1$ be a given constant, $\Lambda$ is defined in (\ref{0.1}), $\eta\leq \va^{\f12+}$,
  then  for any small constant $\sigma>0$, if
 \begin{equation}\Lambda\leq\va^{\f12+\sigma},\label{0.13}\end{equation}
there exists a unique solution $(\mau^\va,\rho^\va)\in \mathcal{B}\times \mathcal{A}$ to the system (\ref{0.3})-(\ref{0.5}) with boundary condition (\ref{0.b1})-(\ref{0.b2}) and the compatibility condition (\ref{0.2}). Moreover, for the remainders $(u,v,\rho)$ defined in the expansion (\ref{ex}), we have the following estimate:
\begin{eqnarray}
\|\mathbf{u}\|_{\mathcal{B}}+\|\rho\|_{\mathcal{A}}\leq \va^{\f\sigma2},\label{main.1}
\end{eqnarray}
and consequently
\begin{eqnarray}
&&|u^\va-u_s|_{L^\infty(\Omega)}+\va^{\f12+}|\nabla(u^\va-u_s)|_{L^\infty}\leq C\va^{\f12+\f\sigma4},\label{main.2}\\
&&|v^\va|_{L^\infty(\Omega)}+\va^{\f12+}|\nabla v^\va|_{L^\infty}\leq C\va^{\f12+\f\sigma4},\label{main.3}\\
&&|\rho^\va-1|_{L^\infty}\leq C\va,
\end{eqnarray}
where  the constant $C$ is independent of  $\va$.
\end{theorem}
The proof of Theorem \ref{th.main} is based on the existence and uniform-in-$\va$ estimate of an approximate system of  the linearized system.  First of all, as there is a second order term $\curl\mau_x$ on the outflow boundary $\{x=L\}$, we can not obtain the existence of solutions to the linearized system directly from the second order elliptic system. Instead, under the boundary condition (\ref{0.b1})-(\ref{0.b2}),  we can first obtain  $\curl\mau$ by taking the curl of the momentum equations. Then  by taking the divergence of the momentum equations, we obtain an elliptic equation for the effective viscous flux $P=a\rho^\gamma-2\dv\mau$. In this step, suitable boundary conditions for $P$ are chosen to ensure that we can back to the original second order equations as well as the compatibility conditions on the corners. Next, based on the the effective viscous flux $P$ and $\curl\mau$ constructed above, we can obtain the density from the mass equation, and the velocity by use of the Helmholtz decomposition of $\mau$.
Here the boundary condition for the potential functions in Helmholtz decomposition are also carefully chosen to ensure the compatibility conditions on the corners in elliptic estimates for the high regularity.

Besides, compared with the incompressible flow, there is a term $\mau^\va\cdot\nabla\rho$ on the left-hand side of the linearized mass equation, leading to a loss of regularity on the right-hand side of the mass equation, it is difficult to show the
existence of a strong solution to the nonlinear system by using the usual fixed point theorem directly. Instead, we will construct a Cauchy consequence $(\mau^n,\rho^n)$ in $(\mathcal{X})^2\times\mathcal{Y}$ to converge to the solution of the nonlinear system, here $\mathcal{X}$ and $\mathcal{Y}$ are subsets of $H^2$ and $ H^1$. For this reason, the estimates for $(\mau,\rho)$ in $(\mathcal{X})^2\times\mathcal{Y}$ must be self enclosed  in the uniform-in-$\va$ estimates of the linearized approximate system.  We achieve this by the delicate analysis on the structure of the linearized system, especially the exploitation of the trace of the solutions to the linear systems on the inflow and outflow part of the boundary.

This paper is organized as follows. In chapter 2, we study the existence and uniform-in-$\va$ estimates of solutions to the linearized system. In section 3, based on the estimate in section 2, we prove the existence of solutions to the nonlinear system as well as the uniform-in-$\va$ estimates of solutions.
%%%%%%%%%%%%%%%%%%%%%%%%%%%%
\renewcommand{\theequation}{\thesection.\arabic{equation}}
\setcounter{equation}{0}
\section{Existence and Uniform-in-$\va$ Estimates of solutions to the linearized system}
To prove Theorem 1.1, we will first consider the linearized system of (\ref{0.3})-(\ref{0.5}), then by taking the limit of a suitable approximating sequence,
we can obtain the existence and uniqueness of  strong solutions to the nonlinear system (\ref{0.3})-(\ref{0.5}) with boundary conditions (\ref{0.b1})-(\ref{0.b2}).
To use the energy method in the following of this section,  we need first to remove the inhomogeneity from the
boundary conditions.

To homogenize the boundary condition, we define the function $\bar{\mathbf{U}}=(\bar u,\bar v)$ satisfying the following systems:
\begin{equation}\label{basic1}\begin{cases}
\Delta^2 \bar u=0\ \text{in}\ \Omega,\\
\bar u|_{y=0}=a_1(0),\ \bar u_y|_{y=0}=b_{0}(x),\\
\bar u|_{y=2}=a_1(2),\  \bar u_y|_{y=2}=b_{1}(x),\\
\bar u|_{x=0}=a_{1}(y),\ \bar u_{xx}|_{x=0}=b_0''(0)y\chi(2y)+(y-2)b_1''(0)\chi(4-2y),\\
\bar u_x|_{x=L}=b_{0}'(L)y\chi(2y)+b_1'(L)(y-2)\chi(4-2y),\\
\bar u_{xxx}|_{x=L}=b_{0}'''(L)y\chi(2y)+b_1'''(L)(y-2)\chi(4-2y).
\end{cases}\end{equation}
and
\begin{equation}\label{basic2}\begin{cases}
\Delta^2\bar v=0\ \text{in}\ \Omega,\\
\bar v|_{y=0,2}=0,\ \bar v_{yy}|_{y=0}=a_2''(0)(x-L)+a''_3(0),\\
 \bar v_{yy}|_{y=2}=a_2''(2)(x-L)+a''_3(2),\\
\bar v_x|_{x=0}=a_{2}(y),\ \bar v_{xxx}|_{x=0}=0,\\
\bar v|_{x=L}=a_{3}(y),\ \bar v_{xx}|_{x=L}=\bar u_{xy}|_{x=L}-a_{4}(y).
\end{cases}\end{equation}
\begin{lemma}
There exists a unique function $\bar{\mau}=(\bar u,\bar v)\in H^4(\Omega)\times H^4(\Omega)$ satisfying system (\ref{basic1}), (\ref{basic2}) and the following estimate:
\begin{equation*}
\|\bar{\mau}\|_{H^4}\lesssim\sum_{i=1}^4\|a_i\|_{H^4}+\|b_1\|_{H^4}+\|b_2\|_{H^4}.
\end{equation*}

\end{lemma}

\begin{proof}

First, Thanks to the compatibility conditions on the corners, if we define
\begin{eqnarray*}h_1(x,y)=&&[a_1(0)+b_0(x)y]\chi(2y)+[a_1(2)-(2-y)b_1(x)]\chi(4-2y),\\
h_2(x,y)=&&\f12y^2[a_2''(0)(x-L)+a''_3(0)]\chi(2y)\\
&&+\f12(y-2)^2[a_2''(2)(x-L)+a''_3(2)]\chi(4-2y),
\end{eqnarray*}
and denote by
\begin{eqnarray*}\hat u=&&\bar u-h_1-[a_1(y)-h_1(0,y)]\chi(\f{4x}L)\triangleq\bar u-u_0,\\
\hat v=&&\bar v-h_2-x[a_2(y)-h_{2x}(0,y)]\chi(\f{4x}L)-[a_3(y)-h_2(L,y)]\\
&&-\f12(x-L)^2\bar v_{xx}(L,y)\chi(\f{4L-4x}{L})
\triangleq\bar u-v_0,
\end{eqnarray*}
then direct computation shows that $(\hat u,\hat v)$ satisfying the following homogeneous boundary value problems:
\begin{equation}\label{iu}\begin{cases}
\Delta^2 \hat u=-\Delta^2 u_0\ \text{in}\ \Omega,\\
\hat u|_{y=0,2}=\hat u_y|_{y=0,2}=\hat u|_{x=0}=\hat u_{xx}|_{x=0}=0,\\
\hat u_x|_{x=L}=\hat u_{xxx}|_{x=L}=0.
\end{cases}\end{equation}
and
\begin{equation}\label{iv}\begin{cases}
\Delta^2 \hat v=-\Delta^2 v_0\ \text{in}\ \Omega,\\
\hat v|_{y=0,2}=\hat v_{yy}|_{y=0,2}=\hat v_x|_{x=0}=\hat v_{xxx}|_{x=0}=0,\\
\hat v|_{x=L}=0,\ \hat v_{xx}|_{x=L}=0,
\end{cases}\end{equation}
where $\Delta^2 u_0,\Delta^2 v_0\in L^2.$

We can prove the existence of solutions to system (\ref{iu}) and (\ref{iv}) by Lax-Milgram Theorem. And under our boundary conditions, the global high regularity in $H^4(\Omega)$ can be proved by odd extension or even extension on the boundary.
\end{proof}
\begin{remark}
If we take $b_0(x)=b_1(x)\equiv0$, i.e. there is no perturbation in  $\curl\mau$ on the boundary $y=0,2$, then the compatibility conditions in (\ref{0.2}) and the process of homogenization above will be simplified.
\end{remark}
It is easy to check that $(\bar u,\bar v)$  satisfy the following boundary conditions:
\begin{equation*}\begin{cases}\bar u_y|_{y=0}=b_0(x),\ \bar u_y|_{y=2}=b_1(x),\ \bar u|_{x=0}=a_1(y),\ \bar v|_{y=0,2}=0,\\
 \bar v|_{x=L}=a_3(y),\ \bar v_x|_{x=0}=a_2(y),\ \partial_x(\bar u_y-\bar v_x)=a_4(y)\end{cases}\end{equation*}

Now we consider the following expansion:
\begin{eqnarray}\label{ex}\begin{cases}
u^\va=u_s+\bar u+\va^{\f12+} u\\
v^\va=\bar v+\va^{\f12+} v\\
\rho^\va=1-\f2\gamma\alpha_2\va\eta^2x+\eta^2\bar \rho+\eta^2\va^{\f12+} \rho.
\end{cases}\end{eqnarray}
here$\eta^2\bar\rho=xh_0(y)$.

Then $(\rho,u,v)$ satisfy the following system:
\begin{eqnarray}
\text{div}\mathbf{u}+\eta^2  u^\va\rho_x+\eta^2  v^\va \rho_y&=&g_0(\rho,u,v)\qquad\text{in}~\Omega,\label{1.1}\\
u_su_x+u_{sy}v-\va\Delta u
-\va\partial_x\mathrm{div}\mathbf{u}+\gamma\rho_x&=&g_1(\rho,u,v)\qquad\text{in}~\Omega,\label{1.2}\\
u_sv_x-\va\Delta v
-\va\partial_y\mathrm{div}\mathbf{u}+\gamma\rho_y&=&g_2(\rho,u,v)\qquad\text{in}~\Omega,\label{1.3}
\end{eqnarray}
with boundary condition
\begin{eqnarray}
\rho_x|_{x=0}=0,\ \rho(0,0)=0,\ u|_{x=0}=v_x|_{x=0}=0,\nonumber\\
v|_{x=L}=\curl\mau_x|_{x=L}=0,\ u_y|_{y=0,2}=v|_{y=0,2}=0.\label{boundary}
\end{eqnarray}
here
\begin{eqnarray}\label{g}\begin{cases}
g_0(\rho,u,v)=g_{01}(\rho,u,v)+g_{0r},\\
g_1(\rho,u,v)=g_{11}(\rho,u,v)+g_{12}(u,v)+g_{1r},\\
g_2(\rho,u,v)=g_{21}(\rho,u,v)+g_{22}(u,v)+g_{2r}\end{cases}
\end{eqnarray}
with
\begin{eqnarray*}
&&g_{01}(\rho,u,v)=\eta^2[\va^{\f12+}\rho\dv\mau+(\f2\gamma\alpha_2\va x-\bar\rho)\dv\mau+\rho\dv\bar\mau+\f2\gamma\alpha_2\va  u-\mau\cdot\nabla\bar\rho]\\
&&g_{11}(\rho,u,v)=(2\alpha_2\va^{\f12-} -\gamma\rho_x)[(\rho^\va)^{\gamma-1}-1]-\va^{-\f12-}\gamma\bar\rho_x(\rho^\va)^{\gamma-1}
-\eta^2\rho(u^\va u^\va_x+v^\va u^\va_y)\\
&&g_{12}(u,v)=- [\va^{\f12+} uu_x+\va^{\f12+} vu_y+\bar uu_x+\bar vu_y+u\bar u_x+v\bar u_y]\\
&&\qquad\qquad-\eta^2(-\f2\gamma\alpha_2\va x+\bar \rho)(u\bar u_x+u^\va u_x+vu^\va_y+\bar vu_y)\\
&&g_{21}(\rho,u,v)=-\gamma\rho_y[(\rho^\va)^{\gamma-1}-1]-\va^{-\f12-}\gamma\bar\rho_y(\rho^\va)^{\gamma-1}-\eta^2 \rho(u^\va v^\va_x+v^\va v^\va_y)\\
&&g_{22}(u,v)=- [\va^{\f12+} uv_x+\va^{\f12+} vv_y+\bar uv_x+\bar vv_y+u\bar v_x+v\bar v_y]\\
&&\qquad\qquad-\eta^2(-\f2\gamma\alpha_2\va x+ \bar \rho)(u^\va v_x+u\bar v_x+v^\va_y v+v_y\bar v)\\
\end{eqnarray*}
and
\begin{eqnarray*}
&&\va^{\f12+} g_{0r}=[\eta^2(\f{2\alpha_2}\gamma\va x- \bar\rho)-1]\dv\bar\mau+\f{2\alpha_2}\gamma\va \eta^2(u_s+\bar u)-\eta^2 [(u_s+\bar u)\bar\rho_x+  \bar v \bar\rho_y]\\
&&\va^{\f12+} g_{1r}=-\bar u\bar u_x-\bar v\bar u_y+\Delta \bar u+\partial_x\mathrm{div}\bar{\mathbf{u}}-(u_s\bar u_x+\bar vu_{sy})\\
&&\qquad\qquad-\eta^2(-\f2\gamma\alpha_2\va x+\bar \rho)[(u_s+\bar u)\bar u_x+\bar v(u_s+\bar u)_y ]\\
&&\va^{\f12+} g_{2r}=\bar u\bar v_x+\bar v\bar v_y+\va\Delta \bar v
+\va\partial_y\mathrm{div}\bar{\mathbf{u}}-\eta^2(-\f2\gamma\alpha_2\va x+ \bar \rho)[(u_s+\bar u)\bar v_x+\bar v\bar v_y)]\\
&&\qquad\qquad-u_s\bar v_x
\end{eqnarray*}

To prove Theorem 1.1, it suffices to prove the existence and uniform-in-$\va$ estimates of a strong solution $(\rho,\mathbf{u})$ to the system
(\ref{1.1})-(\ref{boundary}). To this end, we will construct a sequence that will converge
to the solution of the nonlinear system (\ref{1.1})-(\ref{boundary}).
The solution sequence is defined as follows.
\begin{eqnarray}
\text{div}\mathbf{u}^{n+1}+\eta^2(u_s+\bar u+\va^{\f12+} u^{n})\rho_x^{n+1}+\eta^2(\bar v+\va^{\f12+} v^{n})\rho^{n+1}_y
&=&g_0(\mathbf{u}^{n},\rho^{n}) ,\nonumber\\
u_su_x^{n+1}+u_{sy}v^{n+1}-\va\Delta u^{n+1}
-\va\partial_x\mathrm{div}\mathbf{u}^{n+1}+\gamma\rho_x^{n+1}&=&g_1(\mathbf{u}^{n},\rho^{n}) ,\nonumber\\
u_sv_x^{n+1}-\va\Delta v^{n+1}
-\va\partial_y\mathrm{div}\mathbf{u}^{n+1}+\gamma\rho_y^{n+1}&=&g_2(\mathbf{u}^{n},\rho^{n}) ,\nonumber\\
\rho^{n+1}_x|_{x=0}&=&0 ,\label{3.1}\\
u^{n+1}|_{x=0}=v_x^{n+1}|_{x=0}&=&0,\nonumber\\
v^{n+1}|_{x=L}=\curl\mau^{n+1}_x|_{x=L}&=&0,\nonumber\\
u^{n+1}_y|_{y=0,2}=v^{n+1}|_{y=0,2}&=&0,\nonumber
\end{eqnarray}
where $g_0,g_1,g_2$ are defined in (\ref{g}), and $(\mathbf{u}^0,\rho^0)=(0,0)$.

To show the existence of the solution sequence defined in (\ref{3.1}),  we will first consider the following linear system:
\begin{eqnarray}
\text{div}\mathbf{u}+\eta^2  u^\va\rho_x+\eta^2  v^\va \rho_y&=&g_0\qquad\text{in}~\Omega,\label{linear.1}\\
u_su_x+u_{sy}v-\va\Delta u
-\va\partial_x\mathrm{div}\mathbf{u}+\gamma\rho_x&=&g_1\qquad\text{in}~\Omega,\label{linear.2}\\
u_sv_x-\va\Delta v
-\va\partial_y\mathrm{div}\mathbf{u}+\gamma\rho_y&=&g_2\qquad\text{in}~\Omega,\label{linear.3}
\end{eqnarray}
with boundary condition (\ref{boundary}), where $\mau^\va=(u^\va,v^\va)$, are given functions satisfying
\begin{equation}
\mau^\va\cdot\mathbf{n}=0,\ \text{on}\ y=0,2;\ \|(L-x)\nabla^3\mau^\va\|+\|\mau^\va\|_{W^{2,p}}\leq C,\label{in2}
\end{equation}
and $\ g_0,\ \mathbf{g}=(g_1,g_2)$ are given smooth functions satisfying
\begin{equation}\|\mathbf{g}\|_{H^1}+\|\curl\mathbf{g}\|_{H^1}+|g_2(0,\cdot)|_{H^1}+|g_0(0,\dot)|_{H^2}+\|g_0\|_{H^2}\leq C.\label{in1}\end{equation}

By going to the $\curl$ of (\ref{linear.2})  and (\ref{linear.3}), we have:
\begin{equation*}
u_{sy}\dv\mau+u_s\partial_x(u_y-v_x)+u_{syy}v-\va\Delta(u_y-v_x)=\curl \mathbf{g}, \ \text{in}~\Omega.
\end{equation*}

Then taking the divergence of (\ref{linear.2})  and (\ref{linear.3}), we have:
\begin{equation*}
u_s\partial_x\dv\mau+2u_{sy}v_x+\Delta(\gamma\rho-2\va\dv\mau)=\dv\mathbf{g}\ \text{in}~\Omega.
\end{equation*}
For $\mau=(u,v)$, we define the norm $\|\cdot\|_{\mathcal{X}}$, $\|\cdot\|_{\mathcal{Y}}$  by
\begin{eqnarray}\|\mathbf{u}\|_{\mathcal{X}}:=&&\|\mau\|_{H^1}+\va^{\f12}\|\curl\mau\|_{H^1}+\va\|\mau\|_{H^2}+\va|\dv\mau_y(0,\cdot)|,\label{x}\\
\|\rho\|_{\mathcal{Y}}:=&&\|\rho\|_{H^1}+|\rho_y(0,\cdot)|+\va^{\f12}\eta|\nabla\rho|_{L^\infty_xL^2_y}.\label{y}
\end{eqnarray}

\subsection{Existence of Solutions to an Approximate System}

In our method, to prove the existence of  solutions to the linearized system, we need that $|\nabla\rho(0,\cdot)|_{H^1}$ is well controlled, and consequently we need that $|\curl\mau_x(0,\cdot)|_{H^1}$ is bounded  to construct a compact mapping. This can not be satisfied if we only have $\curl\mau\in W^{2,p}(\Omega)$ in system (\ref{linear.1})-(\ref{linear.3}). To this end,  we will first consider the following approximate system:
\begin{equation}\label{app}\begin{cases}
\text{div}\mau^{app}+\eta^2  u^\va\rho_x^{app}+\eta^2  v^\va \rho_y^{app}=g_0\qquad\text{in}~\Omega, \\
u_su_x^{app}+u_{sy}v^{app}-\va\Delta u^{app}
-\va\partial_x\mathrm{div}\mathbf{u}^{app}+\gamma\rho_x^{app}=g_1\qquad\text{in}~\Omega, \\
u_sv_x^{app}-\va\Delta v^{app}
-\va\partial_y\mathrm{div}\mathbf{u}^{app}+\gamma\rho_y^{app}=g_2\qquad\text{in}~\Omega, \\
2\eta^2\rho_x^{app}(0,y)=\f1{ u^\va}\int_0^y[ (u_y^{app}-v_x^{app})_x^{\alpha}-(u_y^{app}-v_x^{app})_x](0,s)ds,\\
 u^{app}|_{x=0}=v^{app}_x|_{x=0}=v^{app}|_{x=L}=\curl\mau_x^{app}|_{x=L}=0,\\
 u_y^{app}|_{y=0,2}=v^{app}|_{y=0,2}=0,
\end{cases}\end{equation}
here $(u_y^{app}-v_x^{app})_x^{\alpha}=(u_y^{app}-v_x^{app})_x\ast\varrho_\alpha$, $\varrho$ is the standard mollifier, and $\alpha>0$ is a small constant.

\begin{theorem}\label{thapp}
Let $2\leq p<\infty$, for given $\mau^\va,\mathbf{g},g_0$ satisfying (\ref{in2}) and (\ref{in1}), there exists a unique solution $(\mau^{app},\rho^{app})\in (W^{2,p})^2\times H^2$ with the following estimates:
\begin{eqnarray}
&&\|\rho^{app}\|_{\mathcal{Y}}+\va\|\curl\mau^{app}\|_{H^2}+\|\mau^{app}\|_{\mathcal{X}}\nonumber\\
\leq&&C(\|g_0\|_{H^1}+\va|g_{0y}(0,\cdot)|_{L^2}+\|\mathbf{g}\|_{L^2}+\|\curl\mathbf{g}\|_{L^2}+|g_2(0,\cdot)|_{L^2}),\nonumber\\\label{app.2}
\end{eqnarray}
and
\begin{eqnarray}
&&\va^{2-\f2p}\|\mau^{app}\|_{W^{2,p}}+\va^{2-\f2p}\|\curl\mau^{app}\|_{W^{2,p}}]\nonumber\\
\leq&& C\{\|g_0\|_{H^1}+\va|g_{0y}(0,\cdot)|_{L^2}+\|\mathbf{g}\|_{L^2}+\|\curl\mathbf{g}\|_{L^2}+|g_2(0,\cdot)|_{L^2}+\va^{1-\f1p}\eta^{\f2p}|g_2(0,\cdot)|_{L^p}\nonumber\\
&&+\va^{1-\f2p}[\|\mathbf{g}\|_{L^p}+\va|g_{0y}(0,\cdot)|_{L^p}+\va\|g_0\|_{W^{1,p}(\Omega)}+\|\curl\mathbf{g}\|_{L^p}]\},\label{app.3}
\end{eqnarray}
here the constant $C$ is independent of $\alpha,\va,\eta.$
\end{theorem}
\begin{remark}
The only difference between system (\ref{linear.1})-(\ref{linear.3}) with boundary condition (\ref{boundary})  and system (\ref{app}) lies in the boundary condition for $\rho_x|_{x=0}$. In fact, under the new boundary condition, the boundary value for $\rho_y^{app}$ on $x=0$ satisfies
$$\va(u_y^{app}-v_x^{app})_x^{\alpha}-2\va\eta^2(v^\va\rho_y^{app})_y+\gamma\rho_y^{app}=g_2+2\va g_{0y},$$
and by (\ref{in2}) and (\ref{in1}) we can obtain
$$\rho_y^{app}(0,y)\in H^1(0,2),$$
even if we only have $(u_y^{app}-v_x^{app})\in H^2(\Omega)$.
\end{remark}
We will use Leray-Schauder fixed point theory to prove the existence of solutions to system (\ref{app})(for simplicity, we will omit the superscript  in the following of this subsection). For $2<p<\infty$, $\mau=(u,v)$, we denote by
$$\mathfrak{B}=\{\mau\in (W^{1,p}(\Omega))^2,\ \mau_x\in (W^{1,p}(\Omega))^2,\ u|_{x=0}=v|_{y=0,2}=v_x|_{x=0}=v|_{x=L}=0\}.$$
The proof of Theorem \ref{thapp} will be broken into several steps in the following subsections.
\begin{remark}
Here we assume $2<p<\infty$ in $\mathfrak{B}$ as we will need $\nabla\curl\mau\in C(\bar\Omega)$  to ensure the compatibility conditions on the corners. In the following we can prove that $\curl\mau\in W^{2,p}(\Omega)$.
\end{remark}
\subsubsection{Construction of a compact mapping}

For given $\bar\mau=(\bar u,\bar v)\in \mathfrak{B}$, $t\in[0,1]$, we consider the following system:
\begin{equation}\label{app.1}\begin{cases}
\dv\mau+\eta^2u^\va\rho_x+\eta^2v^\va\rho_y=tg_0,\ \text{in}\ \Omega,\\
tu_s\bar u_x+tu_{sy}\bar v-\va\partial_y(u_y-v_x)-2\va\dv\mau_x+\gamma\rho_x=tg_1,\ \text{in}\ \Omega,\\
tu_s\bar v_x+\va\partial_x(u_y-v_x)-2\va\dv\mau_y+\gamma\rho_y=tg_2,\ \text{in}\ \Omega,\\
\rho_x (0,y)=\f1{ 2\eta^2u^\va}\int_0^y[ (u_y -v_x )_x^{\alpha}-(u_y -v_x )_x](0,s)ds,\\
 u |_{x=0}=v _x|_{x=0}=v |_{x=L}=\curl\mau_x |_{x=L}=0,\\
 u_y |_{y=0,2}=v |_{y=0,2}=0.
\end{cases}\end{equation}
In this subsection we will prove the following theorem:
\begin{theorem}\label{com}Let  $\mau^\va,\mathbf{g},g_0$ are given functions satisfying (\ref{in2}) and (\ref{in1}), for given $\bar\mau=(\bar u,\bar v)\in \mathfrak{B}$, $t\in[0,1]$, there exists a compact mapping $T:\mathfrak{B}\times [0,1]\rightarrow\mathfrak{B}$ with $T(\bar\mau,t)=\mau=(u,v)$ be the solution to system (\ref{app.1}), where $\rho$ is defined in (\ref{den.1}).
\end{theorem}
\begin{remark}
As there is a second order term $\partial_x\curl\mau$ on the boundary condition $(\ref{app.1})_5$, we can not obtain the existence directly from the second order system (\ref{app.1}). Instead, under our boundary conditions, we will first prove the existence of $\curl \mau$ and  the effective viscous flux $P=\gamma\rho-2\va\dv\mau$. Here we have chose proper boundary condition for $P$ to make sure that we can back to our original second order system. Then we can obtain the existence of $(\rho,u,v)$ to system (\ref{app.1}) by the mass equation and the Helmholtz decomposition.
\end{remark}

\begin{proof}
We will broke the proof into the following steps.

$\mathbf{Step \ I:}$ Construction of  $H^c=\curl\mau$.

 we consider the following boundary problem for $H^c$:
\begin{equation}\label{curl.2}\begin{cases}
tu_{sy}\dv\bar\mau+tu_s\partial_x(\bar u_y-\bar v_x)+tu_{syy}\bar v-\va\Delta H^c=t\curl\mathbf{g},\\
H^c|_{x=0}=\partial_xH^c|_{x=L}=H^c|_{y=0,2}=0.
\end{cases}\end{equation}
As we have
 $$tu_{sy}\dv\bar\mau+tu_s\partial_x(\bar u_y-\bar v_x)+tu_{syy}\bar v-t\curl\mathbf{g}\in L^p(\Omega),$$
and the compatibility conditions on the corners, we can obtain by the theory of elliptic theories that there exists a unique solution $H^c\in W^{2,p}(\Omega)$ to system (\ref{curl.2}) with the following estimate:
$$\|H^c\|_{W^{2,p}(\Omega)}\leq C(\va)(\|\curl\mathbf{g}\|_{L^p}+\|\bar\mau_x\|_{W^{1,p}(\Omega)}+\|\bar\mau\|_{W^{1,p}(\Omega)}).$$

$\mathbf{Step \ II:}$ Construction of  the effective viscous flux $P=\gamma\rho-2\va\dv\mau$.

 We denote by
\begin{equation}P_0(y)=-2t\va g_0(0,0)+\int_0^y[tg_2-\va H^c_x](0,s)ds,\label{boundary.p}\end{equation}
and define $P$ by the following equation:
\begin{equation}
t(u_s\partial_x\dv\bar\mau+2u_{sy}\bar v_x)-\Delta P=t\dv\mathbf{g},\label{p1}
\end{equation}
with boundary condition:
\begin{equation}P|_{x=0}=P_0;\ P_y|_{y=0,2}=[tg_2-\va H^c_x]|_{y=0,2};\ P_x|_{x=L}=(tg_1-tu_s\bar u_x+\va H^c_y)|_{x=L}.\label{boundary.p1}\end{equation}
Recalling that we already have $\bar u_x\in W^{1,p}(\Omega),\  H^c\in W^{2,p}(\Omega)$, it is easy to check the compatibility conditions on the corners $(0,0)$ and $(0,2)$:
$$P_0'(0)=tg_2(0,0)-\va H^c_x(0,0),\ P_0'(2)=tg_2(0,2)-\va H^c_x(0,2).$$
The theory of elliptic equations implies that there exists a unique solution $P\in H^2(\Omega)$ to equation (\ref{p1}) with boundary condition (\ref{boundary.p1}) and the following estimate:
$$\|P\|_{H^2}\leq C(\va)(\|\dv\mathbf{g}\|_{L^2}+\|\bar\mau_x\|_{H^1(\Omega)}+\|\bar\mau\|_{H^1(\Omega)}).$$

If we denote by $\mathbf{F}=(F_1,F_2)$ with
\begin{eqnarray*}
F_1&=&tu_s\bar u_x+tu_{sy}\bar v-\va\partial_yH^c+\partial_xP-tg_1,\\
F_2&=&tu_s\bar v_x+\va\partial_xH^c+\partial_yP-tg_2,
\end{eqnarray*}
then by (\ref{curl.2}) and (\ref{p1}) we have
$$\dv\mathbf{F}=\curl\mathbf{F}=0,\ \text{in}\ \Omega.$$
Consequently there exists a function $G$ satisfying
$$\nabla G=\mathbf{F},\ \Delta G=0,\ \text{in}\ \Omega.$$
From boundary condition (\ref{boundary.p1}) we have
$$G_y|_{x=0}=G_y|_{y=0,2}=G_x|_{x=L}=0.$$
The theory of elliptic equations imply that
$$G\equiv \text{constant} ,\ \text{in}\ \Omega.$$
That is $F_1=F_2\equiv0$, and we have:
\begin{eqnarray}
tu_s\bar u_x+tu_{sy}\bar v-\va\partial_yH^c+\partial_xP&=&tg_1,\ \text{in}\ \Omega,\label{3.2.3}\\
tu_s\bar v_x+\va\partial_xH^c+\partial_yP&=&tg_2,\ \text{in}\ \Omega.\label{3.3.3}
\end{eqnarray}
$\mathbf{Step \ III:}$  Construction of the density $\rho$ satisfying $(\ref{app.1})_4$.

First for $  0<y<2,$ we denote by $ v^{0}= v^\va(0,y)$, and define a function $\rho^{0,y}(y)$ by the following system:
\begin{eqnarray}
\gamma\rho^{0,y}+2\va\eta^2( v^{0}\rho^{0,y})_{y}=(tg_2+2\va tg_{0y}-\va H_x^{c,\alpha})(0,y),\label{den.y}
\end{eqnarray}
here $H_x^{c,\alpha}=H^c_x\ast\varrho_\alpha$. Since by (\ref{in2}) we have
$$ v^0(0)= v^0(2)=0,\ |v^0|_{H^2}\leq C,$$
for any $2\leq p<\infty$, there exists a unique solution $\rho^{0,y}\in H^1(0,2)$ satisfying equation (\ref{den.y}) and the following estimates:
\begin{eqnarray}|\rho^{0,y}|_{L^p(0,2)}\leq&& C(|g_2(0,\cdot)|_{L^p(0,2)}+\va|g_{0}(0,\cdot)|_{W^{1,p}(0,2)}+\va|H^{c}_x(0,\cdot)|_{L^p})\label{den.y1}
\end{eqnarray}
and
\begin{equation}|\rho_y^{0,y}|_{L^2(0,2)}\leq C(|g_2(0,\cdot)|_{H^1(0,2)}+\va|g_{0}(0,\cdot)|_{H^2(0,2)}+\va|H^{c,\alpha}_{xy}(0,\cdot)|_{L^2(0,2)}),\label{den.yy}\end{equation}
here the constant $C$ is independent of $t,\va,\eta,\alpha$.

Then we consider the following transport equation for $\rho$:
\begin{equation}
\gamma\rho+2\va\eta^2 \mau^\va\cdot\nabla\rho=P+2t\va g_0, \ \rho|_{ x=0}=\int_0^y\rho^{0,y}(s)ds.\label{den.1}
\end{equation}
As $\mau^\va\in W^{2,p}(\Omega),\ P+2t\va g_0\in H^2(\Omega),\ \mau^\va\cdot\mathbf{n}=0$ on $y=0,2$, and
 $$\|\mau^\va-\mathbf{U}_0\|_{w^{2,p}}\ll1,$$
we obtain that  there exists a unique solution $\rho\in H^2(\Omega)$  to problem (\ref{den.1}) satisfying:
\begin{eqnarray*}
&&\|\rho\|_{H^2(\Omega)}+|\nabla^2\rho|_{L^\infty_xL^2_y}\leq C(\va,\eta)(\|P\|_{H^2}+\|g_0\|_{H^2}+|\rho_y^{0,y}|_{H^1(0,2)})\nonumber\\
\leq&&C(\va,\eta,\alpha)(\|\mathbf{g}\|_{H^1(\Omega)}+\|g_0\|_{H^2}+\|\bar\mau_x\|_{H^1(\Omega)}+|g_2(0,\cdot)|_{H^1(0,2)}+\va|g_{0}(0,\cdot)|_{H^2(0,2)})\nonumber\\
\label{den.app}
\end{eqnarray*}

Now we differentiate (\ref{den.1}) with respect to $y$ to have
\begin{equation}\gamma\rho_y+2\va\eta^2( u^\va\rho_x+ v^\va\rho_y)_y=P_y+2t\va g_{0y},\label{den.y2}\end{equation}
which Combined with the boundary condition (\ref{boundary.p1}) and (\ref{den.y}) implies
$$2\va\eta^2(  u^\va\rho_x)_y=\va( H_x^{c,\alpha}-H^c_x),\ \text{on}\ x=0.$$
Besides, as $v^\va(0,0)=0$, we have by (\ref{boundary.p}) and (\ref{den.1})
$$\rho(0,0)=(  v^\va\rho_y)(0,0)=(P+2t\va g_0)(0,0)=0,$$
then from equation (\ref{den.1}), there holds
$$  u^\va\rho_x(0,0)=0.$$
Consequently we have
\begin{equation}2\eta^2\rho_x (0,y)=\f1{ u^\va}\int_0^y[ H^{c,\alpha}_x-H^c_x](0,s)ds\label{boundary.rhox}\end{equation}
and
$$|\rho_x(0,\cdot)|_{H^1}\leq C\eta^{-2}|(H_x^{c,\alpha}-H_x^c)(0,\cdot)|_{L^2(\Omega)}\lesssim\f\alpha{\eta^2}|H^c_x(0,\cdot)|_{L^2(\Omega)}.$$

$\mathbf{Step \ IV:}$ Construction of $\mau=(u,v)$

We write $\mau$ in the form of Helmholtz decomposition:
\begin{equation}u=\psi_y+\phi_x,\ v=-\psi_x+\phi_y.\label{hel}\end{equation}
Here $\phi(x,y)$ is defined by
\begin{equation}\label{div.2}\begin{cases}
2\va\Delta\phi=-P+\gamma\rho,\ \text{in}\ \Omega,\\
\phi_x|_{x=0}=\phi_y|_{y=0,2}=\phi|_{x=L}=0,
\end{cases}\end{equation}
and $\psi(x,y)$ is defined by
\begin{equation}\label{curl.1}\begin{cases}
\Delta\psi=H^c,\ \text{in}\ \Omega,\\
\psi|_{x=0}=\psi|_{y=0,2}=0,\\
\psi_x|_{x=L}=0.
\end{cases}\end{equation}
First it is easy to check that
\begin{equation*}
u|_{x=0}=0,\ v|_{y=0,2}=0,\ v|_{x=L}=\partial_x\curl\mau|_{x=L}=0, v_x|_{x=0}=(-H^c+u_y)|_{x=0}=0.\label{boundary.1}
\end{equation*}
Next, by Lemma \ref{lemdiv} and Lemma \ref{lemcurl} in the following of this subsection, we have
$$\mau\in W^{2,p}(\Omega), \ \mau_x\in  W^{1+\f1p,p}(\Omega).$$
As $W^{1+\f1p,p}\subset\subset  W^{1,p},$ we have proved that $T$ is a compact mapping from $\mathfrak{B}\times[0,1]$ to $\mathfrak{B}$.

$\mathbf{Step \ IV:}$ Back to the original system.

Combing (\ref{3.2.3}),(\ref{3.3.3}), (\ref{den.1}), (\ref{boundary.rhox}) and $(\ref{div.2})_1$ we obtain that the triple $(u,v,\rho)$ solves the system (\ref{app.1}), and thus we finish the proof of Theorem \ref{com}.

\end{proof}

\begin{lemma}\label{lemdiv}
There exists a unique solution $\phi$ to system (\ref{div.2}) with the following estimate:
$$\|\phi\|_{W^{3,p}(\Omega)}+\|\phi_{xx}\|_{W^{1+\f1p,p}(\Omega)}\leq C(\|P\|_{H^2}+\|\rho\|_{H^1}+\|\rho_x\|_{H^1}+|\rho(L,\cdot)|_{H^1(0,2)}),$$
here $C$ is a constant depends only on $\Omega$.
\end{lemma}
\begin{remark}
In Lemma \ref{lemdiv}, as we do not have the compatibility conditions on the corners $(L,0)$ and $(L,2)$ for the second order derivatives of $\phi$, we can not obtain the global regularity for $\phi$ in $H^4(\Omega)$ even through we already  have  $-P+\gamma\rho\in H^2(\Omega)$.
\end{remark}
\begin{proof}
First of all, as $-P+\gamma\rho\in H^2(\Omega)\hookrightarrow W^{1,p}(\Omega)$, for any $2\leq p<\infty$, and  the compatibility conditions at the corners:
$$\phi_{y}(L,0)=\phi_{y}(L,2)=0, $$
there exists a unique solution $\phi\in W^{2,p}(\Omega)$ to system (\ref{div.2}) satisfying
$$\va\|\phi\|_{W^{2,p}(\Omega)}\leq \|P-\gamma\rho\|_{L^p(\Omega)}.$$
Then we differentiate equation $(\ref{div.2})_1$ with respect to $y$ to have
\begin{equation}\label{div.3}
2\va\Delta\phi_y=-P_y+\gamma\rho_y,\ \text{in}\ \Omega,\end{equation}
with boundary condition
\begin{equation*}
\phi_{xy}|_{x=0}=\phi_y|_{y=0,2}=\phi_y|_{x=L}=0.
\end{equation*}
It is easy to check the compatibility condition at the corners:
$$\phi_{xy}(0,0)=\phi_{xy}(0,2)= \phi_{y}(L,0)=\phi_{y}(L,2)=0. $$
So we have $\phi_y\in W^{2,p}(\Omega)$ with
$$\ \va\|\phi_y\|_{W^{2,p}}\leq \|P_y-\gamma\rho_y\|_{L^p(\Omega)},$$
 which combined with equation $(\ref{div.2})_1$ implies $\phi\in W^{3,p}(\Omega)$ and
 $$\ \va\|\phi\|_{W^{3,p}}\leq \|P\|_{W^{1,p}}+\|\rho\|_{W^{1,p}(\Omega)}.$$

Next we differentiate equation $(\ref{div.2})_1$ with respect to $x$ twice to have
\begin{equation}\label{div.4}
2\va\Delta\phi_{xx}=-P_{xx}+\gamma\rho_{xx},\ \text{in}\ \Omega,\end{equation}
with boundary condition
$$\begin{cases}\phi_{xxy}|_{y=0,2}=0,\\
2\va \phi_{xxx}|_{x=0}=[(-P+\gamma\rho)_x-2\va\phi_{xyy}]|_{x=0}=(-P+\gamma\rho)_x|_{x=0},\\
2\va\phi_{xx}|_{x=L}=[(-P+\gamma\rho)-2\va\phi_{yy}]|_{x=L}=(-P+\gamma\rho)|_{x=L}.
\end{cases}$$
As  $-P_x+\gamma\rho_x\in H^1(\Omega)$, for any subset $V\subset\bar\Omega$ with $V\cap\{(L,y),0\leq y\leq2\}=\emptyset$, we have $\phi_{xx}\in H^2(V)$. To obtain the global estimate, we first take the even extension of $\phi_{xx},P_{xx},\rho_{xx}$ with respect to $y=0$, and denote the even extension of a function $\vartheta$ by
$$\vartheta^*_{xx}(x,y)=\begin{cases}\vartheta_{xx}(x,y), \ 0<y<2,\\\vartheta_{xx}(x,-y), -2\leq y\leq0,\end{cases}$$
and  $\Omega^*=(0,L)\times(-2,2)$. Then we have
$$\Delta\phi^*=-P_{xx}^*+\gamma\rho_{xx}^*,\ \text{in}\ \Omega^*.$$
Thanks to the homogeneous Neumann boundary condition on $y=0$,  we have $\phi_{xx}^*\in H^2(V^*)$ for any subset $V^*\subset\bar\Omega^*$ with $V^*\cap\{(L,y),-2\leq y\leq2\}=\emptyset$.

On the other hand, as we have $P\in H^2(\Omega)$ and $\rho_{yy}\in L^\infty_xL^2_y$, the Sobolev imbedding theory implies   $(-P+\rho)_y(L,y)\in L^p(0,2)$, for any $2\leq p<\infty$. The even tension of $(-P+\rho)(L,y)$,  denoted by $(-P^*+\rho^*)(L,y)$ is in $W^{1,p}(-2,2)$. By the standard estimate for the Dirichlet problem of Laplace Operator in smooth domain(see e.g. \cite{JK})  we obtain
$$\phi^*_{xx}\in W^{1+\f1p,p}(\Omega_0^*), \ \Omega_0^*=(0,L)\times(-2+\sigma,2-\sigma),$$
 Consequently we have
$$\phi_{xx}\in W^{1+\f1p,p}(\Omega_0), \ \Omega_0=(0,L)\times(0,2-\sigma).$$

Similarly we can take the even extension of $\phi_{xx}$ with respect to $y=2$ to obtain
$$\phi_{xx}\in W^{1+\f1p,p}(\Omega_1), \ \Omega_1=(0,L)\times(\sigma,2).$$
Combing above, we have proved that
$$\phi_{xx}\in W^{1+\f1p,p}(\Omega),$$
and we complete the proof.
\end{proof}

\begin{lemma}\label{lemcurl}
There exists a unique solution $\psi\in H^4(\Omega)$ to system (\ref{curl.1}) with estimate:
\begin{equation*}\|\psi\|_{W^{4,p}(\Omega)}\leq C\|H^c\|_{W^{2,p}(\Omega)}.\end{equation*}
\end{lemma}
\begin{proof}
First, as we have the compatibility conditions on the corners, we have
$$\psi\in W^{2,p}(\Omega)$$
Next, we take derivative of equation $(\ref{curl.1})_1$ with respect to $x$ to have
\begin{equation}\Delta\psi_x=H^c_x,\ \text{in}\ \Omega,\ \psi_x|_{x=L}=\psi_x|_{y=0,2}=0.\label{curl.3}\end{equation}
besides, as
$$\Delta\psi|_{x=0}=H^c|_{x=0}=0, \psi|_{x=0}=0,$$
we have
$$\psi_{xx}|_{x=0}=0.$$
by the fact that $H^c\in W^{2,p}$ and the compatibility conditions on the corners, we have
$$\psi_x\in W^{2,p}(\Omega).$$
Then back to equation $(\ref{curl.1})_1$, we have
$$\psi\in W^{3,p}(\Omega).$$
Next we take derivative of equation (\ref{curl.3}) with respect to $x$  to have
$$\Delta\psi_{xx}=H^c_{xx},\ \text{in}\ \Omega,\ \psi_{xx}|_{x=0}=\psi_{xx}|_{y=0,2}=0.$$
besides, as
$$\Delta\psi_x|_{x=L}=H^c_x|_{x=L}=0, \psi_x|_{x=L}=0,$$
we have
$$\psi_{xxx}|_{x=L}=0.$$
by the fact that $H^c\in W^{2,p}$ and the compatibility conditions on the corners, we have
$$\psi_{xx}\in W^{2,p}(\Omega).$$
 Taking derivative of equation (\ref{curl.3}) with respect to $y$, we have
 $$\psi_{xyyy}=-\psi_{xxxy}+H^c_{xy}\in L^p(\Omega),$$
 Finally,  by taking derivative of equation $(\ref{curl.1})_1$ with respect to $y$ twice we have
 $$\psi_{yyyy}=-\psi_{xxyy}+H^c_{yy}\in L^p(\Omega).$$
 Combing above, we have proved $\psi\in W^{4,p}(\Omega)$ with
 $$\|\psi\|_{W^{4,p}(\Omega)}\leq \|H^c\|_{W^{2,p}},$$
 and we complete the proof.
\end{proof}

\subsubsection{Uniform-in-t estimates in Hilbert space}

To use the Leray-Schauder fixed point theory,  we still need to prove the uniform-in-$t$ estimate of the fixed point satisfying $T(\mau,t)=\mau$, i.e.:
\begin{eqnarray}
\dv\mau+\eta^2u^\va\rho_x+\eta^2v^\va\rho_y&=&tg_0,\ \text{in}\ \Omega,\label{3.3}\\
tu_s u_x+tu_{sy} v-\va\partial_y(u_y-v_x)-2\va\partial_x\dv\mau+\gamma\partial_x\rho&=&tg_1,\ \text{in}\ \Omega,\label{3.4}\\
tu_s v_x+\va\partial_x(u_y-v_x)-2\va\partial_y\dv\mau+\gamma\partial_y\rho&=&tg_2,\ \text{in}\ \Omega,\label{3.5}
\end{eqnarray}
with boundary condition:
\begin{equation}\label{boundary.2}\begin{cases}
\rho_x (0,y)=\f1{ 2\eta^2u^\va}\int_0^y[ (u_y -v_x )_x^{\alpha}-(u_y -v_x )_x](0,s)ds,\\
 u |_{x=0}=v _x|_{x=0}=v |_{x=L}=\curl\mau_x |_{x=L}=0,\\
 u_y |_{y=0,2}=v |_{y=0,2}=0.
 \end{cases}\end{equation}
 For simplicity we first introduce the following notations:
\begin{eqnarray*}
A_{1,t}&=&t\|u_x\|+t\| v_{x}\|,\label{A1}\\
A_{2,t}&=&t\|(u_y-v_x)_x\|+t^{\f12}\va^{\f12}|(u_y-v_x)_x(0,\cdot)|+t^{\f12}\va^{\f12}|(u_y-v_x)_y(L,\cdot)|.\label{A2}
\end{eqnarray*}
As there is a term $\mau^\va\cdot\nabla\rho$
on the left-hand side of the mass equation (\ref{3.1})$_1$, so that there will be loss of regularity on the right-hand side of the equation. To prove the existence of solutions to  the nonlinear system in next section, we will have to find the fixed point in the space $(\rho,\mau)\in \mathcal{Y}\times (\mathcal{X})^2$, which does not involve the $H^2$ norm of $\rho$. So in the uniform-in-$\va$ estimates for the linearized system, the estimates for $(\rho,\mau)$ in the space $\mathcal{Y}\times (\mathcal{X})^2$ should be self enclosed. More precisely, we will prove the following theorem in this subsection:
\begin{theorem}\label{th3.1} For given $\mau^\va, g_0,\mathbf{g}=(g_1,g_2)$ satisfying (\ref{in2}) and (\ref{in1}), $t\in[0,1]$, let $(u,v,\rho)$ be solutions to system (\ref{3.3})-(\ref{boundary.2}), $A_i,i=1,2$ are defined above, then
we have:
\begin{eqnarray*}
&&A_{1,t}+A_{2,t}+\|\rho\|_{\mathcal{Y}}+\va\|\curl\mau\|_{H^2}+t^{\f12}\va^{\f12}\|\curl\mau\|_{H^1}+t\|\mau\|_{H^1}+\va\|\mau\|_{H^2}\nonumber\\
&&+\va|\dv\mau_y(0,\cdot)|\nonumber\\
\leq&&C(\|g_0\|_{H^1}+\va|g_{0y}(0,\cdot)|_{L^2}+\|\mathbf{g}\|_{L^2}+\|\curl\mathbf{g}\|_{L^2}+|g_2(0,\cdot)|_{L^2}),
\end{eqnarray*}
here the constant $C$ is independent of $t,\alpha,\va,\eta.$
\end{theorem}
The proof of Theorem \ref{th3.1} will be broken into several steps. First by going to the curl of (\ref{3.4}) and (\ref{3.5}) we have
\begin{equation}\label{curl.4}\begin{cases}
tu_{sy}\dv\mau+tu_s\partial_x( u_y- v_x)+tu_{syy} v-\va\Delta (u_y-v_x)=t\curl\mathbf{g},\\
H^c|_{x=0}=\partial_xH^c|_{x=L}=H^c|_{y=0,2}=0.
\end{cases}\end{equation}
Then combing (\ref{3.3})-(\ref{3.5}) we have
\begin{eqnarray}
\gamma\rho_x+2\va\eta^2(u^\va\rho_{xx}+u^\va_x\rho_x+v^\va\rho_{xy}+v^\va_x\rho_y)&&=\va(u_y-v_x)_y\nonumber\\
&& +t(g_1+2\va g_{0x}-u_su_x-u_{sy}v),\label{rhox}\\
\gamma\rho_y+2\va\eta^2(u^\va\rho_{xy}+v^\va\rho_{yy}+u^\va_y\rho_x+v^\va_y\rho_y)=&&-\va(u_y-v_x)_x+t(g_2+2\va g_{0y}-u_sv_x).\nonumber\\\label{rhoy}
\end{eqnarray}
We have the following estimates for the density:
\begin{lemma}\label{den}
Let $2\leq p<\infty$,\  $(u,v,\rho)$ be solutions to system (\ref{3.3})-(\ref{boundary.2}), and the assumptions in Theorem \ref{th3.1} are satisfied, then we have:
\begin{eqnarray*}
\|\rho_x\|_{L^p}+\|\rho_y\|_{L^p}&&+|\rho_y(0,\cdot)|_{L^p}+\va^{\f1p}\eta^{\f2p}|\nabla\rho|_{L^\infty_xL_y^p}\leq C [\va^{\f1p}\eta^{\f2p}|g_{2}(0,\cdot)|_{L^p(0,2)}\nonumber\\
&&+\|\mathbf{g}\|_{L^p}+\va|g_{0y}(0,\cdot)|_{L^p}+\va\|H^c\|_{H^2}+t\|\mau_x\|_{L^p(\Omega)}+\va\|g_0\|_{W^{1,p}(\Omega)}]
\end{eqnarray*}
here the constant $C$ is independent of $t,\va,\eta,\alpha.$
\end{lemma}
\begin{proof}

First for any $x_0\in[0,L]$, we multiply (\ref{rhoy}) with $\rho_y|\rho_y|^{p-2}$ to have
\begin{eqnarray*}
&&\int_0^2\int_0^{x_0}\gamma\rho_y^pdxdy+2\va\eta^2\int_0^2\int_0^{x_0}\rho_y|\rho_y|^{p-2}(  u^\va\rho_{xy}+ u_y^\va\rho_x+  v^\va\rho_{yy}+  v^\va_y\rho_y)dxdy\nonumber\\
=&&\int_0^2\int_0^{x_0}[(1+(2-\f{2}p)\eta^2\va   v^\va_y-\f2p\va\eta^2  u^\va_x)\rho_y^p+2\va\eta^2 u_y^\va\rho_x\rho_y|\rho_y|^{p-2}]dxdy\nonumber\\
&&+\f2p\va\eta^2\int_0^2  u^\va|\rho_y|^p|_{x=x_0}dy-\f2p\va\eta^2\int_0^2  u^\va|\rho_y|^p|_{x=0}dy\nonumber\\
=&&\int_0^2\int_0^{x_0}\rho_y|\rho_y|^{p-2}[-\va(u_y-v_x)_x+t(g_2+2\va g_{0y}-u_sv_x)]dxdy.
\end{eqnarray*}
By (\ref{den.y1}) we have
$$|\rho_y(0,\cdot)|_{L^p}\lesssim |g_2(0,\cdot)|_{L^p}+\va|g_{0y}(0,\cdot)|_{L^p}+\va|H^c_x(0,\cdot)|_{L^p},$$
and consequently
\begin{eqnarray}
&&\|\rho_y\|_{L^p(\Omega)}+\va^{\f1p}\eta^{\f2p}\sup_{x\in[0,L]}|\rho_y(x,\cdot)|_{L^p_y}\nonumber\\
\lesssim&& \va\eta^2\|\rho_x\|_{L^p(\Omega)}+\va^{\f1p}\eta^{\f2p}|\rho_y(0,\cdot)|_{L^p}+\va\|\curl\mau\|_{H^2(\Omega)}+\va\|g_0\|_{W^{1,p}(\Omega)}+t\|v_x\|_{L^p(\Omega)}\nonumber\\
&&+\|g_2\|_{L^p(\Omega)}\nonumber\\
\lesssim&&\va\eta^2\|\rho_x\|_{L^p(\Omega)}+\va^{\f1p}\eta^{\f2p}[|g_2(0,\cdot)|_{L^p(0,2)}+\va|g_{0y}(0,\cdot)|_{L^p}]+\va\|\curl\mau\|_{H^2(\Omega)}\nonumber\\
&&+\va\|g_0\|_{W^{1,p}(\Omega)}+t\|v_x\|_{L^p(\Omega)}+\|g_2\|_{L^p(\Omega)}\label{den.y3}
\end{eqnarray}
Similarly   we multiply (\ref{rhox}) with $\rho_x|\rho_x|^{p-2}$ to obtain
\begin{eqnarray*}
&&\int_0^2\int_0^{x_0}\gamma\rho_x^pdxdy+2\va\eta^2\int_0^2\int_0^{x_0}[( (1-\f1p) u^\va_x- \f1pv^\va_y)\rho_x^p+2v_x^\va\rho_x\rho_y|\rho_x|^{p-2}]dxdy\nonumber\\
&&+\f2p\va\eta^2\int_0^2  u^\va\rho_x^p|_{x=x_0}dy\nonumber\\
=&&\f2p\va\eta^2\int_0^2  u^\va\rho_x^p|_{x=0}dy+t\int_0^2\int_0^{x_0}\rho_x|\rho_x|^{p-2}(g_1+2\va g_{0x}-u_su_x-u_{sy}v)dxdy\nonumber\\
&&+\va\int_0^2\int_0^{x_0}\rho_x|\rho_x|^{p-2}(u_y-v_x)_ydxdy.
\end{eqnarray*}
Using (\ref{boundary.rhox}) we have
\begin{eqnarray}
&&\|\rho_x\|_{L^p}+\va^{\f1p}\eta^{\f2p}\sup_{x\in[0,L]}|\rho_x(x,\cdot)|_{L^p_y}\nonumber\\
\lesssim&& \va^{\f1p}\eta^{\f2p}|\rho_x(0,\cdot)|_{L^p}+\va\|(u_y-v_x)_y\|_{L^p}+\va\|g_0\|_{W^{1,p}}+\va^2\eta^2\|\rho_y\|_{L^p}+t\|\mau_x\|_{L^p}\nonumber\\
&&+\|g_1\|_{L^p}\nonumber\\
\lesssim &&\va\|(u_y-v_x)\|_{H^2}+\va\|g_0\|_{W^{1,p}}+\va^2\eta^2\|\rho_y\|_{L^p}+t\|\mau_x\|_{L^p}+\|g_1\|_{L^p}.\label{den.x3}
\end{eqnarray}
Combing (\ref{den.y3}) and (\ref{den.x3}), we have proved Lemma \ref{den}.

\end{proof}

\begin{lemma}\label{lema2}
Let $(u,v,\rho)$ be solutions to system (\ref{3.3})-(\ref{boundary.2}), and the assumptions in Theorem \ref{th3.1} are satisfied, then we have:
\begin{equation*}
\va\|u_y-v_x\|_{H^2}+A_{2,t}\leq C( LA_{1,t}+\|g_0\|_{L^2}+\|\curl\mathbf{g}\|_{L^2}+\eta^2\|\nabla\rho\|_{L^2}),
\end{equation*}
here the constant $C$ is independent of $t,\alpha,\va,\eta.$
\end{lemma}
\begin{proof}
First we multiply (\ref{curl.4}) with $t(u_y-v_x)_x$ to have
\begin{eqnarray*}
&&t\int_\Omega[tu_s\partial_x( u_y- v_x)+tu_{syy} v-\va(u_y-v_x)_{xx}-\va(u_y-v_x)_{yy}](u_y-v_x)_xdxdy\nonumber\\
=&&t^2\int_\Omega [u_s(u_y-v_x)_x^2+u_{syy}v(u_y-v_x)_x]dxdy+\f12t\va\int_0^2(u_y-v_x)_x^2|_{x=0}dy\nonumber\\
&&+\f12t\va\int_0^2(u_y-v_x)_y^2|_{x=L}dy\nonumber\\
=&&t\int_\Omega(t\curl\mathbf{g}-tu_{sy}\dv\mau)(u_y-v_x)_xdxdy
\end{eqnarray*}
Consequently we have
\begin{eqnarray}
&&A_{2,t}=t\|(u_y-v_x)_x\|+t^{\f12}\va^{\f12}|(u_y-v_x)_x(0,\cdot)|+t^{\f12}\va^{\f12}|(u_y-v_x)_y(L,\cdot)|\nonumber\\
\lesssim&& t(\|v\|+\|\curl\mathbf{g}\|_{L^2}+\|\dv\mau\|_{L^2})\nonumber\\
\lesssim&& LA_{1,t}+\|\curl\mathbf{g}\|_{L^2}+\eta^2\|\nabla\rho\|_{L^2}+\|g_0\|_{L^2}.\label{curl.7}
\end{eqnarray}

Then we multiply (\ref{curl.4}) with $\va(u_y-v_x)_{xx}$ to have
\begin{eqnarray*}
&&\va^2\int_\Omega(u_y-v_x)_{xx}\Delta(u_y-v_x)dxdy=\va^2\int_\Omega[(u_y-v_x)_{xx}^2+(u_y-v_x)_{xy}^2]dxdy\\
=&&\va\int_\Omega(u_y-v_x)_{xx}[t\curl\mathbf{g}-tu_{sy}\dv\mau-tu_s\partial_x( u_y- v_x)-tu_{syy} v]dxdy\\
\lesssim &&t\va\|(u_y-v_x)_{xx}\|(\|v\|+\|\curl\mathbf{g}\|_{L^2}+\|\dv\mau\|+\|\partial_x(u_y-v_x)\|)\nonumber\\
\lesssim&&\va\|(u_y-v_x)_{xx}\|( LA_{1,t}+\|\curl\mathbf{g}\|_{L^2}+t\|\partial_x(u_y-v_x)\|+\eta^2\|\nabla\rho\|+\|g_0\|),
\end{eqnarray*}
which combined with equation (\ref{curl.4}) and (\ref{curl.7}) implies
$$\va\|u_y-v_x\|_{H^2}\lesssim LA_{1,t}+\|g_0\|_{L^2}+\|\curl\mathbf{g}\|_{L^2}+\eta^2\|\nabla\rho\|_{L^2}.$$
Combing above, we have proved Lemma \ref{lema2}.
\end{proof}

\begin{lemma}\label{lema1}
Let $(u,v,\rho)$ be solutions to system (\ref{3.3})-(\ref{boundary.2}), and the assumptions in Theorem \ref{th3.1} are satisfied, then we have:
\begin{eqnarray*}
A_{1,t}^2\leq&& C(LA_{2,t}^2+L\va^2\|\nabla \mau_x\|^2+\|g_0\|_{H^1}^2+|g_2(0,\cdot)|_{L^2}^2+\va|g_{0}(0,\cdot)|_{H^1(0,2)}\\
&&+\eta^2\|\nabla\rho\|_{L^2}^2+\eta^4|\rho_x(L,\cdot)|^2+\|\mathbf{g}\|_{L^2}^2+\|\rho_x\|\|g_0\|_{L^2}),
\end{eqnarray*}
here the constant $C$ is independent of $t,\alpha,\va,\eta.$
\end{lemma}
\begin{proof}
We multiply (\ref{3.4}) with $tu_x$ and (\ref{3.5}) with $tv_x$ to have
\begin{eqnarray}
&&t\int_\Omega[tu_s u_x+tu_{sy} v-\va\partial_y(u_y-v_x)-2\va\partial_x\dv\mau+\partial_x\rho]u_xdxdy\nonumber\\
&&+t\int_\Omega[tu_s v_x+\va\partial_x(u_y-v_x)-2\va\partial_y\dv\mau+\partial_y\rho]v_xdxdy\nonumber\\
=&&t^2\int_\Omega [u_s(u_x^2+v_x^2)+u_{sy}u_xv]dxdy+t\va\int_\Omega(u_yu_{xy}-2u_{xx}u_x-v_{xy}u_x)dxdy\nonumber\\
&&+t\va\int_\Omega(-v_{xx}v_x-u_{xy}v_x+2v_yv_{xy})dxdy-t\int_0^2\rho_yv|_{x=0}dy+t\int_\Omega\rho_x\dv\mau dxdy\nonumber\\
=&&t^2\int_\Omega [u_s(u_x^2+v_x^2)+u_{sy}u_xv]dxdy+\f12t\va\int_0^2u_y^2|_{x=L}dy+t\int_\Omega\rho_x\dv\mau dxdy\nonumber\\
&&-t\va\int_0^2u_x^2|_{x=L}dy+t\va\int_0^2u_x^2|_{x=0}dy-\f12t\va\int_0^2v_x^2|_{x=L}dy-t\va\int_0^2v_y^2|_{x=0}dy\nonumber\\
&&-t\int_0^2\rho_yv|_{x=0}dy\nonumber\\
=&&t^2\int_\Omega [u_s(u_x^2+v_x^2)+u_{sy}u_xv]dxdy+\f12t\va\int_0^2(u_y-v_x)(u_y+v_x)|_{x=L}dy\nonumber\\
&&+t\int_\Omega\rho_x\dv\mau dxdy+t\va\int_0^2\dv\mau(u_x-v_y)|_{x=0}dy-t\va\int_0^2u_x^2|_{x=L}dy\nonumber\\
&&-t\int_0^2\rho_yv|_{x=0}dy\nonumber\\
=&&t^2\int_\Omega(g_1u_x+g_2v_x)dxdy,
\end{eqnarray}
Consequently  we have
\begin{eqnarray}
&&t^2\int_\Omega u_s(u_x^2+v_x^2)dxdy\nonumber\\
\lesssim&&t^2\|u_x\|\|v\|+t\|\rho_x\|\|\dv\mau\|+t|\rho_y(0,\cdot)||v(0,\cdot)|+t\va |u_{x}(L,\cdot)|^2+A_{1,t}\|\mathbf{g}\|\nonumber\\
&&+t\va|(u_y-v_x)(L,\cdot)||(u_y+v_x)(L,\cdot)|+t\va|\dv\mau(0,\cdot)||(u_x-v_y)(0,\cdot)|
\end{eqnarray}
here by equation (\ref{3.3}) we have
\begin{eqnarray}
&&|\dv\mau(0,\cdot)|+|u_x(L,\cdot)|\nonumber\\
\lesssim&& t|g_0(0,\cdot)|+\eta^2|\nabla\rho(0,\cdot)|+ t|g_0(L,\cdot)|+\eta^2|\rho_x(L,\cdot)|,
\end{eqnarray}
besides,
\begin{eqnarray*}&&t\va|(u_y-v_x)(L,\cdot)||(u_y+v_x)(L,\cdot)|+t\va|\dv\mau(0,\cdot)||(u_x-v_y)(0,\cdot)|\\
\lesssim&&tL\va\|(u_y-v_x)_x\|\|\nabla\mau_x\|+t\va|\dv\mau(0,\cdot)|(|u_x(L,\cdot)|+L^{\f12}\|\nabla\mau_x\|)\\
\lesssim&&L^{\f12}\va^2\|\nabla\mau_x\|^2+LA_{2,t}^2+\eta^2|g_2(0,\cdot)|^2+\|g_0\|_{H^1}^2+\eta^4|\rho_x(L,\cdot)|^2+\eta^2\va^2|H^c_x(0,\cdot)|^2
\end{eqnarray*}
and
\begin{eqnarray}
&&t^2\|u_x\|\|v\|+t\|\rho_x\|\|\dv\mau\|+t|\rho_y(0,\cdot)||v(0,\cdot)|\nonumber\\
\lesssim&& LA_{1,t}^2+\|\rho_x\|\|g_0\|_{L^2}+\eta^2\|\nabla\rho\|^2+L A_{1,t}(|g_2(0,\cdot)|_{L^2(0,2)}+\va|g_{0}(0,\cdot)|_{H^1(0,2)}\nonumber\\
&&+\va|H^{c}_x(0,\cdot)|_{L^2})
\end{eqnarray}
Combing above, using Lemma \ref{den} and Lemma \ref{lema2} we have
\begin{eqnarray*}
A_{1,t}^2\lesssim&& LA_{2,t}^2+L\va^2\|\nabla \mau_x\|^2+\|g_0\|_{H^1}^2+|g_2(0,\cdot)|_{L^2}^2+\va|g_{0}(0,\cdot)|_{H^1(0,2)}\\
&&+\eta^2\|\nabla\rho\|_{L^2}^2+\eta^4|\rho_x(L,\cdot)|^2+\|\mathbf{g}\|_{L^2}^2+\|\rho_x\|\|g_0\|_{L^2},
\end{eqnarray*}
and we complete the proof.
\end{proof}

$\mathbf{Proof\ of\ Theorem\ \ref{th3.1}}:$

First by  Lemma \ref{den}-\ref{lema1}  we have
\begin{eqnarray*}
&&A_{1,t}+A_{2,t}+\va\|H^c\|_{H^2}\\
\lesssim&&L^{\f12}\va\|\nabla \mau_x\|+\|g_0\|_{H^1}+|g_2(0,\cdot)|_{L^2}+\va|g_{0y}(0,\cdot)|_{L^2}+\|\mathbf{g}\|_{L^2}+\|\curl\mathbf{g}\|_{L^2}\\
&&+\eta^2\|\nabla\rho\|_{L^2}+\eta^2|\rho_x(L,\cdot)|+\|\rho_x\|^{\f12}\|g_0\|_{L^2}^{\f12}
\end{eqnarray*}
Taking $p=2$ in Lemma \ref{den}, using the fact that $\eta<\va^{\f12+}$ we obtain
\begin{eqnarray}
&&A_{1,t}+A_{2,t}+|\rho_y(0,\cdot)|_{L^p}+\va\|u_y-v_x\|_{H^2}+\|\nabla\rho\|_{L^2}+\va^{\f12}\eta|\nabla\rho|_{L^\infty_xL_y^2}\nonumber\\
\lesssim&&L^{\f12}\va\|\nabla u_x\|+\|g_0\|_{H^1}+\va|g_{0y}(0,\cdot)|_{L^2}+\|\mathbf{g}\|_{L^2}+\|\curl\mathbf{g}\|_{L^2}+|g_2(0,\cdot)|_{L^2}.\nonumber\\
\label{2.1}
\end{eqnarray}
Using  the momentum equations (\ref{3.4})-(\ref{3.5}) we have
\begin{eqnarray}
\va\|\nabla\dv\mau\|\lesssim \|\nabla\rho\|+\va\|u_y-v_x\|_{H^1}+t\|\mau_x\|+\|\mathbf{g}\|,
\end{eqnarray}
Now if we write equation (\ref{3.5}) in the form
\begin{equation}\label{sysv}\begin{cases}
-\va\Delta v=tg_2-\gamma\rho_y+\va\dv\mau_y-tu_sv_x,\ \text{in}\ \Omega,\\
v|_{x=L}=v|_{y=0,2}=v_x|_{x=L}=0,
\end{cases}\end{equation}
then thanks to the compatibility conditions on the corners we have by the theory of elliptic equations
\begin{equation}\va\|v\|_{H^2}\lesssim \|g_2\|_{L^2}+\|\rho_y\|_{L^2}+\va\|\dv\mau_y\|_{L^2}+t\|v_x\|_{L^2}\label{v.3}.\end{equation}
Recalling  that $u_y=(u_y-v_x)+v_x$, we obtain
$$\va\|u_y\|_{H^1}\leq \va\|v\|_{H^2}+\va\|u_y-v_x\|_{H^1},$$
which combined with equation (\ref{3.4}),(\ref{2.1})-(\ref{v.3}) implies
\begin{eqnarray}
\va\|\mau\|_{H^2}\lesssim \|g_0\|_{H^1}+\va|g_{0y}(0,\cdot)|_{L^2}+\|\mathbf{g}\|_{L^2}+\|\curl\mathbf{g}\|_{L^2}+|g_2(0,\cdot)|_{L^2}.
\label{u1}
\end{eqnarray}
besides, we have
$$t\|u_y\|\leq t\|u_y-v_x\|+t\|v_x\|\leq t\|(u_y-v_x)_x\|+t\|v_x\|\leq A_{1,t}+A_{2,t}, $$
and by equation (\ref{3.3})
$$t\|v_y\|\lesssim t\|u_x\|+\eta^2\|\nabla\rho\|+\|g_0\|\leq A_{1,t}+\eta^2\|\nabla\rho\|+\|g_0\|.$$
Moreover, by (\ref{boundary.p1})  we have
\begin{eqnarray*}
&&\va|\dv\mau_y(0,\cdot)|=|(P_y-\gamma\rho_y)(0,\cdot)|\\
\lesssim&&|g_2(0,\cdot)|+\va|(u_y-v_x)_x(0,\cdot)|+|\rho_y(0,\cdot)|\lesssim |g_2(0,\cdot)|+\va\|u_y-v_x\|_{H^2}+\va|g_{0y}(0,\cdot)|.
\end{eqnarray*}
Finally we multiply (\ref{curl.4}) with $t(u_y-v_x)$ to have
\begin{eqnarray*}
&&t\int_\Omega[u_s\partial_x(u_y-v_x)-\va\Delta(u_y-v_x)](u_y-v_x)dxdy\nonumber\\
=&&\f12t\int_0^2u_s(u_y-v_x)^2|_{x=L}dy+t\va\int_\Omega|\nabla(u_y-v_x)|^2dxdy\nonumber\\
=&&t^2\int_\Omega(\curl\mathbf{g}-u_{sy}\dv\mau-u_{syy}v)(u_y-v_x)dxdy\nonumber\\
\lesssim&&L^{\f12}t^2(\|\curl\mathbf{g}\|+\|\dv\mau\|+\|v\|)\|(u_y-v_x)_x\|\leq L^{\f12}A_{2,t}(\|\curl\mathbf{g}\|+t\|\mau\|_{H^1})
\end{eqnarray*}
Combing above, using the fact that $L\ll1$ we obtain
\begin{eqnarray*}
&&A_{1,t}+A_{2,t}+\va\|\curl\mau\|_{H^2}+\|\nabla\rho\|_{L^2}+\va^{\f12}\eta|\nabla\rho|_{L^\infty_xL_y^2}+\va|\dv\mau_y(0,\cdot)|\nonumber\\
&&+\va\|\mau\|_{H^2}+t\|\mau\|_{H^1}+t^{\f12}\va^{\f12}\|u_y-v_x\|_{H^1}\nonumber\\
\lesssim&&\|g_0\|_{H^1}+\va|g_{0y}(0,\cdot)|_{L^2}+\|\mathbf{g}\|_{L^2}+\|\curl\mathbf{g}\|_{L^2}+|g_2(0,\cdot)|_{L^2},
\end{eqnarray*}
and we have finish the proof of Theorem \ref{th3.1}.

\subsubsection{Proof of Theorem \ref{thapp}}

To prove Theorem \ref{thapp}, we still need the uniform-in-$t$ estimate for $\mau$ in $W^{2,p}(\Omega)$.

First by Theorem \ref{3.1} and Sobolev imbedding theorem we have
\begin{eqnarray}
&&t^{\f2p}\va^{1-\f2p}\|\mau_x\|_{L^p}+t^{\f2p}\va^{1-\f2p}\|(u_y-v_x)_x\|_{L^p}+t^{\f2p}\va^{1-\f2p}\|\dv\mau\|_{L^p}\nonumber\\
\leq&&(t\|\mau_x\|)^{\f2p}(\va\|\mau_x\|_{H^1})^{1-\f2p}+ (\va\|u_y-v_x\|_{H^2})^{1-\f2p}(t\|(u_y-v_x)_x\|_{L^2})^{\f2p}\nonumber\\
&&+(t\|\mau\|_{H^1})^{\f2p}(\va\|\mau\|_{H^2})^{1-\f2p}\nonumber\\
\lesssim &&\|g_0\|_{H^1}+\va|g_{0y}(0,\cdot)|_{L^2}+\|\mathbf{g}\|_{L^2}+\|\curl\mathbf{g}\|_{L^2}+|g_2(0,\cdot)|_{L^2},\label{p.5}
\end{eqnarray}
Using the momentum equations (\ref{3.4})-(\ref{3.5}) again we have for any $2\leq p<\infty$,
\begin{eqnarray*}
\va\|\nabla\dv\mau\|_{L^p}\lesssim \|\nabla\rho\|_{L^p}+\va\|u_y-v_x\|_{H^2}+t\|\mau_x\|_{L^p}+\|\mathbf{g}\|_{L^p},
\end{eqnarray*}
which combined with (\ref{p.5}) and Lemma \ref{den} implies
\begin{eqnarray}
&&\va^{1-\f2p}[\va\|\nabla\dv\mau\|_{L^p}+\|\rho_x\|_{L^p}+\|\rho_y\|_{L^p}+|\rho_y(0,\cdot)|_{L^p}+\va^{\f1p}\eta^{\f2p}|\nabla\rho|_{L^\infty_xL_y^p}]\nonumber\\
\lesssim&&\|g_0\|_{H^1}+\va|g_{0y}(0,\cdot)|_{L^2}+\|\mathbf{g}\|_{L^2}+\|\curl\mathbf{g}\|_{L^2}+|g_2(0,\cdot)|_{L^2}+\va^{1-\f1p}\eta^{\f2p}|g_2(0,\cdot)|_{L^p}\nonumber\\
&&+\va^{1-\f2p}[\|\mathbf{g}\|_{L^p}+\va|g_{0y}(0,\cdot)|_{L^p}+\va\|g_0\|_{W^{1,p}(\Omega)}].\label{p.1}
\end{eqnarray}
By system (\ref{sysv}) and the theory of elliptic equations
\begin{eqnarray}
&&\va^{2-\f2p}\|v\|_{W^{2,p}}\nonumber\\
\lesssim &&\va^{1-\f2p}(\|g_2\|_{L^p}+\|\rho_y\|_{L^p}+\va\|\dv\mau_y\|_{L^p}+t\|v_x\|_{L^p})\nonumber\\
\lesssim&& \|g_0\|_{H^1}+\va|g_{0y}(0,\cdot)|_{L^2}+\|\mathbf{g}\|_{L^2}+\|\curl\mathbf{g}\|_{L^2}+|g_2(0,\cdot)|_{L^2}+\va^{1-\f1p}\eta^{\f2p}|g_2(0,\cdot)|_{L^p}\nonumber\\
&&+\va^{1-\f2p}[\|\mathbf{g}\|_{L^p}+\va|g_{0y}(0,\cdot)|_{L^p}+\va\|g_0\|_{W^{1,p}(\Omega)}].\label{p.4}
\end{eqnarray}
Since
\begin{eqnarray*} \va^{2-\f2p}\|u_y\|_{W^{1,p}}\leq\va^{2-\f2p}(\|u_y-v_x\|_{W^{1,p}}+\|v\|_{2,p})
\leq\va\|u_y-v_x\|_{H^2}+\va^{2-\f2p}\|v\|_{2,p},
\end{eqnarray*}
combining equation (\ref{3.4}) we have
\begin{eqnarray}
&&\va^{2-\f2p}\|u\|_{2,p}\nonumber\\
\lesssim&&\|g_0\|_{H^1}+\va|g_{0y}(0,\cdot)|_{L^2}+\|\mathbf{g}\|_{L^2}+\|\curl\mathbf{g}\|_{L^2}+|g_2(0,\cdot)|_{L^2}+\va^{1-\f1p}\eta^{\f2p}|g_2(0,\cdot)|_{L^p}\nonumber\\
&&+\va^{1-\f2p}[\|\mathbf{g}\|_{L^p}+\va|g_{0y}(0,\cdot)|_{L^p}+\va\|g_0\|_{W^{1,p}(\Omega)}].\label{p.6}
\end{eqnarray}
Finally by system (\ref{curl.4}) and (\ref{p.5})  we have
\begin{eqnarray}
&&\va^{2-\f2p}\|u_y-v_x\|_{W^{2,p}}\nonumber\\
\leq&& \va^{1-\f2p}\|tu_{sy}\dv\mau+tu_s\partial_x( u_y- v_x)+tu_{syy} v-t\curl\mathbf{g}\|_{L^p}\nonumber\\
\lesssim&& \|g_0\|_{H^1}+\va|g_{0y}(0,\cdot)|_{L^2}+\|\mathbf{g}\|_{L^2}+\|\curl\mathbf{g}\|_{L^2}+|g_2(0,\cdot)|_{L^2}+\va^{1-\f1p}\eta^{\f2p}|g_2(0,\cdot)|_{L^p}\nonumber\\
&&+\va^{1-\f2p}[\|\mathbf{g}\|_{L^p}+\va|g_{0y}(0,\cdot)|_{L^p}+\va\|g_0\|_{W^{1,p}(\Omega)}+\|\curl\mathbf{g}\|_{L^p}].\label{pc1}
\end{eqnarray}

Now we are ready to prove Theorem \ref{thapp} by use of Leray-Schauder fixed point theory.
By Theorem \ref{com}, there exists a compact  mapping $T$ from $\mathfrak{B}\times[0,1]\rightarrow\mathfrak{B}$  with $T(\bar\mau,t)=\mau$ be the solution to system (\ref{app.1}).

For any $\bar\mau\in\mathfrak{B}$, when $t=0$, we have by (\ref{curl.2}) and (\ref{p1}) that
 $$H^c=P\equiv0.$$
Then by (\ref{den.y}) and (\ref{den.1}) we have
$$\rho^{0,y}=\rho\equiv0.$$
Consequently by (\ref{hel})-(\ref{curl.1}) we have $T(\bar\mau,0)=\mau\equiv0.$

 Finally from  (\ref{p.4}) and (\ref{p.6}), there exists a constant $C$ independent of $t$ such that
  $$\|\mau\|_{W^{2,p}}\leq C(\va)(\|\curl\mathbf{g}\|_{L^p}+\|g_0\|_{H^1}+\|g_0\|_{H^2}+\|\mathbf{g}\|_{H^1})$$
  for all $(\mau,t)\in\mathfrak{B}\times [0,1]$ satisfying $\mau=T(\mau,t)$. By Leray-Schauder fixed point theorem (c.f.\cite{G-T}), the mapping $T_1(\bar\mau)=T(\bar\mau,1)$ of $\mathfrak{B}$ into itself has a fixed point with estimates in Theorem \ref{th3.1} and (\ref{p.4})-(\ref{pc1}), thus we have complete the proof of Theorem \ref{thapp}.

\subsection{Higher regularity for the approximate solutions}

In this subsection we will prove some higher regularity for the approximate solutions.   For simplicity, we will still omit the superscript  in the following of this subsection.  First have the following boundary estimates for  $\curl\mau$.
\begin{lemma}\label{lem.h1} For given $\mau^\va, g_0,\mathbf{g}=(g_1,g_2)$ satisfying (\ref{in2}) and (\ref{in1}),
let $(u,v,\rho)$ be solutions to the approximate  system (\ref{app}), then we have:
\begin{eqnarray*}
&&\va^{\f32}|(u_y-v_x)_{xx}(L,\cdot)|+\va^{\f32}|(u_y-v_x)_{xy}(0,\cdot)|\nonumber\\
\leq &&C(\va|g_{0y}(0,\cdot)|_{L^2}+\|\mathbf{g}\|_{L^2}+|g_2(0,\cdot)|_{L^2} +\|g_0\|_{H^1}+\|\curl\mathbf{g}\|+\va^{\f12}\|\curl\mathbf{g}\|_{H^1})\end{eqnarray*}
here the constant $C$ is independent of $\va,\eta$.
\end{lemma}
\begin{proof}
First we differentiate $(\ref{curl.4})_1$ with respect to $x$, taking $t=1$ to have
\begin{equation}
u_{sy}\dv\mau_x+u_s\partial_{xx}( u_y- v_x)+u_{syy} v_x-\va\Delta (u_y-v_x)_x=\curl\mathbf{g}_x.\label{curl.5}
\end{equation}
Then we multiply (\ref{curl.5}) with $\va^2(u_y-v_x)_{xx}$ to obtain
\begin{eqnarray}
&&-\va^3\int_\Omega[(u_y-v_x)_{xxx}+(u_y-v_x)_{xyy}](u_y-v_x)_{xx}dxdy\nonumber\\
=&&-\f12\va^3\int_0^2(u_y-v_x)_{xx}^2|_{x=L}dy+\f12\va^3\int_0^2(u_y-v_x)_{xx}^2|_{x=0}dy\nonumber\\
&&-\f12\va^3\int_0^2(u_y-v_x)_{xy}^2|_{x=0}dy\nonumber\\
=&&\va^2\int_\Omega[\curl\mathbf{g}_x-u_{sy}\dv\mau_x-u_s\partial_{xx}( u_y- v_x)-u_{syy} v_x](u_y-v_x)_{xx}dxdy,\nonumber
\end{eqnarray}
which combined with Theorem \ref{th3.1} and equation $(\ref{curl.4})_1$ implies
\begin{eqnarray*}
&&\va^3\int_0^2(u_y-v_x)_{xx}^2|_{x=L}dy+\va^3\int_0^2(u_y-v_x)_{xy}^2|_{x=0}dy\nonumber\\
\lesssim&& \va |(u_{sy}\dv\mau+u_s\partial_x( u_y- v_x)+u_{syy} v-\curl\mathbf{g})(0,\cdot)|_{L^2}^2+\va^2\|(u_y-v_x)_{xx}\|_{L^2}^2\nonumber\\
&&+\va^2(\|v_x\|_{L^2}^2+\|\dv\mau_x\|_{L^2}^2+\|\curl\mathbf{g}\|_{H^1}^2)\nonumber\\
\lesssim&& \va^2|g_{0y}(0,\cdot)|_{L^2}^2+\|\mathbf{g}\|_{L^2}^2+|g_2(0,\cdot)|_{L^2}^2 +\|g_0\|_{H^1}^2+\va\|\curl\mathbf{g}\|_{H^1}^2,
\end{eqnarray*}
then Lemma \ref{lem.h1} follows immediately.
\end{proof}
\begin{lemma}\label{u.3}Let $(u,v,\rho)$ be solutions to the approximate system (\ref{app}), and the assumptions in Lemma \ref{lem.h1} are satisfied,  then we have:
\begin{eqnarray*}|u_x(0,\cdot)|+\va^{\f12}\|\mau\|_{H^2}&&\leq C\{\va^{\f12}\|g_0\|_{H^2}+\|\mathbf{g}\|_{L^2}+|\mathbf{g}(0,\cdot)|_{L^2} \\ &&+\|g_0\|_{H^1}+\|\curl\mathbf{g}\|+L^{\f12}\va^{\f12}\eta^2|\rho_{xy}(L,\cdot)|)
+ \va^{\f12}\|\mathbf{g}\|_{H^1}\}
\end{eqnarray*}
here the constant $C$ is independent of $\alpha,\va,\eta$.
\end{lemma}
\begin{proof}
First we differentiate $(\ref{app})_1$ with respect to $y$ to have
$$\dv\mau_y+\eta^2u_{sy}\rho_x+\eta^2u_s\rho_{xy}+\eta^2(v^\va\rho_y)_y=g_{0y}.$$
As $v(L,y)=v^\va(L,y)=0$, $\eta<\va^{\f12+}$, we have by Theorem \ref{th3.1}
\begin{eqnarray}
\va^{\f12}|u_{xy}(L,\cdot)|=&&\va^{\f12}|(g_{0y}-\eta^2u_{sy}\rho_x-\eta^2u_s\rho_{xy})(L,\cdot)|\nonumber\\
\lesssim &&\va^{\f12}\|g_0\|_{H^2}+ \va^{\f12}\|\mathbf{g}\|_{H^1}+\va^{\f12}\eta^2|\rho_{xy}(L,\cdot)|),\label{uxy}
\end{eqnarray}
Then we differentiate $(\ref{app})_3$ with respect to $x$ to have
\begin{equation}
u_sv_{xx}-\va\Delta v_x-\va\partial_y\dv\mau_x+\partial_y\rho_x=g_{2x}.\label{3.5.1}
\end{equation}

We multiply $(\ref{app})_2$ with $u_{xx}$ and (\ref{3.5.1}) with $v_{x}$ to have
\begin{eqnarray}
&&\int_\Omega[u_s u_x+u_{sy} v-\va u_{yy}-2\va u_{xx}-\va v_{xy}+\partial_x\rho]u_{xx}dxdy\nonumber\\
&&-\int_\Omega[u_s v_{xx}+\va\partial_{xx}(u_y-v_x)-2\va\partial_{xy}\dv\mau+\partial_{xy}\rho]v_{x}dxdy\nonumber\\
=&&\f12\int_0^2u_su_x^2|_{x=L}dy-\f12\int_0^2u_su_x^2|_{x=0}dy-\int_0^2u_{sy}u_xv|_{x=0}dy\nonumber\\
&&-\int_\Omega u_{sy}u_xv_xdxdy+\va\int_0^2u_yu_{xy}|_{x=L}dy-\va\int_\Omega(u_{xy}^2+2u_{xx}^2+u_{xx}v_{xy})dxdy\nonumber\\
&&-\f12\int_0^2v_x^2|_{x=L}dy+\va\int_0^2u_{xy}v_x|_{x=L}dy-\va\int(v_{xx}^2+u_{xx}v_{xy}+2v_{xy}^2)dxdy\nonumber\\
&&+\int_\Omega\rho_x\partial_x\dv\mau dxdy\nonumber\\
=&&-\f12\int_0^2u_su_x^2|_{x=0}dy-\f12\int_0^2v_x^2|_{x=L}dy+\f12\int_0^2u_su_x^2|_{x=L}dy\nonumber\\
&&-\int_0^2u_{sy}u_xv|_{x=0}dy-\int_\Omega u_{sy}u_xv_xdxdy+\va\int_0^2(u_y+v_x)u_{xy}|_{x=L}dy\nonumber\\
&&-\va\int_\Omega[u_{xy}^2+u_{xx}^2+v_{xx}^2+v_{xy}^2+(\partial_x\dv\mau)^2]dxdy+\int_\Omega\rho_x\partial_x\dv\mau dxdy\nonumber\\
=&&\int_\Omega(g_1u_{xx}+g_{2x}v_x)dxdy
\end{eqnarray}
here
\begin{eqnarray*}
&&\int_\Omega\rho_x\partial_x\dv\mau dxdy=\int_\Omega\rho_x\partial_x(g_0-\eta^2u^\va\rho_x-\eta^2v^\va\rho_y)dxdy\\
=-&&\f12\eta^2\int_0^2u^\va\rho_x^2|_{x=L}dy+\f12\eta^2\int_0^2u^\va\rho_x^2|_{x=0}dy+\int_\Omega g_{0x}\rho_xdxdy\\
&&-\f12\eta^2\int_\Omega[(u_x^\va-v^\va_y)\rho_x^2+v^\va_x\rho_x\rho_y]dxdy.
\end{eqnarray*}
Combining above, using Theorem \ref{3.1} we have
\begin{eqnarray*}
&&|u_x(0,\cdot)|^2+|v_x(L,\cdot)|^2+\eta^2|\rho_x(L,\cdot)|^2+\va\|\nabla\mau_x\|^2\nonumber\\
\lesssim&&|v(0,\cdot)|^2+\|\mau_x\|^2+\va|(u_y+v_x)(L,\cdot)||u_{xy}(L,\cdot)|+\|\rho_x\|\|g_0\|_{H^1}+|g_1(0,\cdot)|\nonumber\\
&&+\eta^2|\rho_x(0,\cdot)^2|+\eta^2\|\nabla\rho\|^2\nonumber\\
\lesssim&&\va|g_{0y}(0,\cdot)|_{L^2}+\|\mathbf{g}\|_{L^2}+|\mathbf{g}(0,\cdot)|_{L^2} +\|g_0\|_{H^1}+\|\curl\mathbf{g}\|\nonumber\\
&&+L^{\f12}\va(\|u_{xy}\|+\|v_{xx}\|)|u_{xy}(L,\cdot)|.
\end{eqnarray*}
Consequently from   (\ref{uxy}) we have
\begin{eqnarray*}
&&|u_x(0,\cdot)|+\va^{\f12}\|\mau_x\|_{H^1}\\
\lesssim &&\va^{\f12}\|g_0\|_{H^2}+\|\mathbf{g}\|_{L^2}+|\mathbf{g}(0,\cdot)|_{L^2} +\|g_0\|_{H^1}+\|\curl\mathbf{g}\|+L^{\f12}\va^{\f12}\eta^2|\rho_{xy}(L,\cdot)|)\\
&&+ \va^{\f12}\|\mathbf{g}\|_{H^1}.
\end{eqnarray*}
Finally taking $t=1$ in Theorem \ref{3.1}, we obtain
\begin{eqnarray*}&&\va^{\f12}\|u_{yy}\|\leq \va^{\f12}(\|(u_y-v_x)_y\|+\|v_{xy}\|)\nonumber\\
\lesssim&&\va^{\f12}\|g_0\|_{H^2}+\|\mathbf{g}\|_{L^2}+|\mathbf{g}(0,\cdot)|_{L^2} +\|g_0\|_{H^1}+\|\curl\mathbf{g}\|+L^{\f12}\va^{\f12}\eta^2|\rho_{xy}(L,\cdot)|)\\
&&+ \va^{\f12}\|\mathbf{g}\|_{H^1}.
\end{eqnarray*}
Combing above, we have finish the proof of the Lemma.
\end{proof}

Based on Lemma \ref{lem.h1}, we can obtain the following estimate for the density and the effective viscous flux:
\begin{lemma}\label{den.lem}
Let $(u,v,\rho)$ be solutions to the approximate system (\ref{app}), and the assumptions in Lemma \ref{lem.h1} are satisfied,  then we have:
\begin{eqnarray*}
&&\va\eta\|\rho_{xx}\|+\va^{\f12}\|\nabla\rho_{y}\|+\va\eta^2|\rho_{xx}|_{L^\infty_xL^2_y}+\va\eta|\nabla\rho_{y}|_{L^\infty_xL^2_y}+\va^{\f12}|\rho_{yy}(0,\cdot)|\\
\leq&&C\{\va^{\f12}\|g_{0}\|_{H^2}+\|\mathbf{g}\|_{L^2}+|g_2(0,\cdot)|_{L^2} +\|g_0\|_{H^1}+\|\curl\mathbf{g}\|+\va^{\f12}\|\curl\mathbf{g}\|_{H^1}\nonumber\\
&&+\va|g_{0yy}(0,\cdot)|+\va^{\f12}\|\mathbf{g}\|_{H^1}+\va^{\f12}|g_{2y}(0,\cdot)|\},
\end{eqnarray*}
here the constant $C$ is independent of $\va,\eta,\alpha.$
\end{lemma}
\begin{proof}
First by  (\ref{den.y1}), (\ref{boundary.rhox}), (\ref{rhox}), Theorem \ref{3.1} and Lemma \ref{u.3} we find
\begin{eqnarray*}
&&\va\eta^2|\rho_{xx}(0,\cdot)|\nonumber\\
\leq&& |\rho_x(0,\cdot)|_{H^1(0,2)}+\va^2\eta^2|\rho_y(0,\cdot)|+\va|g_{0x}(0,\cdot)|+|g_1(0,\cdot)|+|u_x(0,\cdot)|+|v(0,\cdot)|\nonumber\\
\lesssim&&  \va^{\f12}\|g_0\|_{H^2}+\|\mathbf{g}\|_{L^2}+|\mathbf{g}(0,\cdot)|_{L^2}+\|g_0\|_{H^1}+\|\curl\mathbf{g}\|+L^{\f12}\va^{\f12}\eta^2|\rho_{xy}(L,\cdot)|\nonumber\\
&&+\va^{\f12}\|\mathbf{g}\|_{H^1}.
\end{eqnarray*}
Then we differentiate (\ref{rhox}) with respect to $x$ to have
\begin{eqnarray}&&\gamma\rho_{xx}+2\va\eta^2(  u^\va\rho_{xx}+  u^\va_x\rho_x+  v^\va_x\rho_y+\  v^\va\rho_{xy})_x\nonumber\\
=&&\va(u_y-v_x)_{xy} +(g_1+2\va g_{0x}-u_su_x-u_{sy}v)_x.\label{den.xx}\end{eqnarray}

For any fixed $x_0\in[0,L]$,  we multiply (\ref{den.xx}) with  $\va\eta^2\rho_{xx}$ to have
\begin{eqnarray*}
&&\va\eta^2\int_0^2\int_0^{x_0}\rho_{xx}^2dxdy+\va^2\eta^4\int_0^2  u^\va\rho_{xx}|_{x=x_0}dy-\va^2\eta^4\int_0^2  u^\va\rho_{xx}|_{x=0}dy\nonumber\\
&&+\va^2\eta^4\int_0^2\int_0^{x_0}[(3  u^\va_x-  v_y^\va)\rho_{xx}^2+2(  u^\va_{xx}\rho_x+  v^\va_{xx}\rho_y+2  v_x^\va\rho_{xy})\rho_{xx}]dxdy\nonumber\\
=&&\va\eta^2\int_0^2\int_0^{x_0}\rho_{xx}[\va(u_y-v_x)_{xy} +(g_1+2\va g_{0x}-u_su_x-u_{sy}v)_x]dxdy.
\end{eqnarray*}
Consequently we have
\begin{eqnarray}
&&\va^{\f12}\eta\|\rho_{xx}\|+\va\eta^2\sup_{x\in[0,L]}|\rho_{xx}(x,\cdot)|_{L^2_y}\nonumber\\
\lesssim&& \va\eta^2|\rho_{xx}(0,\cdot)|+\eta^3\va^2(\|\rho_x\|+\|\rho_y\|+\|\rho_{xy}\|)+\va^{\f32}\eta\|u_y-v_x\|_{H^2}\nonumber\\
&&+\va^{\f12}\eta(\va\|g_{0}\|_{H^2}+\|\mathbf{g}\|_{H^1}+\|\mau_{x}\|_{H^1})\nonumber\\
\lesssim&&\va^{\f12}\|g_0\|_{H^2}+\|\mathbf{g}\|_{L^2}+|\mathbf{g}(0,\cdot)|_{L^2}+\|g_0\|_{H^1}+\|\curl\mathbf{g}\|+L^{\f12}\va^{\f12}\eta^2|\rho_{xy}(L,\cdot)|\nonumber\\
&&+\va^{\f12}\|\mathbf{g}\|_{H^1}.\label{den.2}
\end{eqnarray}
Similarly we  differentiate (\ref{rhoy}) with respect to $x$ and multiply $\va\rho_{xy}$ to have
\begin{eqnarray*}
&&\va\int_0^2\int_0^{x_0}\rho_{xy}^2dxdy+\va^2\eta^2\int_0^2 u^\va\rho_{xy}^2|_{x=x_0}dy-\va^2\eta^2\int_0^2 u^\va\rho_{xy}^2|_{x=0}dy\nonumber\\
&&+\va^2\eta^2\int_0^2\int_0^{x_0}[(u^\va_x+v^\va_y)\rho_{xy}^2+2(  u_y^\va\rho_{xx}+  u^\va_{xy}\rho_x+  v_{xy}^\va\rho_y+  v_y^\va\rho_{yy})\rho_{xy}]dxdy\nonumber\\
=&&\va\int_0^2\int_0^{x_0}\rho_{xy}[-\va(u_y-v_x)_{xx}+g_{2x}+2\va g_{0xy}-u_sv_{xx}]dxdy,
\end{eqnarray*}
which implies
\begin{eqnarray*}
&&\va^{\f12}\|\rho_{xy}\|+\va\eta\sup_{x\in[0,L]}|\rho_{xy}(x,\cdot)|_{L^2_y}\nonumber\\
\lesssim &&\va\eta|\rho_{xy}(0,\cdot)|+\va\eta^2(\|\nabla\rho\|+\|\rho_{xx}\|+\|\rho_{yy}\|)+\va^{\f32}\|u_y-v_x\|_{H^2}+\va^{\f12}\|v_{xx}\|\nonumber\\
&&+\|g_{2x}\|+\va\|g_{0xy}\|\nonumber\\
\lesssim&&\va^{\f12}\|g_0\|_{H^2}+\|\mathbf{g}\|_{L^2}+|\mathbf{g}(0,\cdot)|_{L^2}+\|g_0\|_{H^1}+\|\curl\mathbf{g}\|+L^{\f12}\va^{\f12}\eta^2|\rho_{xy}(L,\cdot)|\nonumber\\
&&+\va^{\f12}\|\mathbf{g}\|_{H^1}+\eta^2\va\|\rho_{yy}\|.
\end{eqnarray*}
Recalling that $\eta<\va^{\f12+}$, we have
\begin{eqnarray}
&&\va^{\f12}\|\rho_{xy}\|+\va\eta\sup_{x\in[0,L]}|\rho_{xy}(x,\cdot)|_{L^2_y}\nonumber\\
\lesssim&&\va^{\f12}\|g_0\|_{H^2}+\|\mathbf{g}\|_{L^2}+|\mathbf{g}(0,\cdot)|_{L^2}+\|g_0\|_{H^1}+\|\curl\mathbf{g}\|+\va^{\f12}\|\mathbf{g}\|_{H^1}
+\eta^2\va\|\rho_{yy}\|\nonumber\\
\label{den.3}
\end{eqnarray}
Next we  differentiate (\ref{rhoy}) with respect to $y$ and multiply $\va\rho_{yy}$ to have
\begin{eqnarray*}
&&\va\int_0^2\int_0^{x_0}\rho_{yy}^2dxdy+\va^2\eta^2\int_0^2  u^\va\rho_{yy}^2|_{x=x_0}dy-\va^2\eta^2\int_0^2  u^\va\rho_{yy}^2|_{x=0}dy\nonumber\\
&&+\va^2\eta^2\int_0^2\int_0^{x_0}[(3  v^\va_y-  u^\va_x)\rho_{yy}^2+2(2  u_y^\va\rho_{xy}+  u_{yy}^\va\rho_x+ v^\va_{yy}\rho_y)\rho_{yy}dxdy\nonumber\\
=&&\va\int_0^2\int_0^{x_0}\rho_{yy}[-\va(u_y-v_x)_{xy}+g_{2y}+2\va g_{0yy}-(u_sv_x)_y]dxdy.
\end{eqnarray*}
By equation (\ref{den.yy}) and Lemma \ref{lem.h1} we have
\begin{eqnarray*}
&&\va^{\f12}|\rho_{yy}(0,\cdot)|\lesssim\va^{\f12}(|g_2(0,\cdot)|_{H^1}+\va|g_0(0,\cdot)|_{H^2}+\va|H_{xy}^c(0,\cdot)|)\\
\lesssim&&\va^{\f12}\|g_{0}\|_{H^2}+\|\mathbf{g}\|_{L^2}+|g_2(0,\cdot)|_{L^2} +\|g_0\|_{H^1}+\|\curl\mathbf{g}\|+\va^{\f12}\|\curl\mathbf{g}\|_{H^1}\\
&&+\va|g_{0yy}(0,\cdot)|+\va^{\f12}|g_{2y}(0,\cdot)|.
\end{eqnarray*}
Consequently we have
\begin{eqnarray}
&&\va^{\f12}\|\rho_{yy}\|+\va\eta\sup_{x\in[0,L]}|\rho_{yy}(x,\cdot)|_{L^2_y}\nonumber\\
\lesssim &&\va\eta|\rho_{yy}(0,\cdot)|+\va\eta^2(\|\nabla\rho\|+\|\rho_{xy}\|)+\va^{\f32}\|u_y-v_x\|_{H^2}+\va\|g_0\|_{H^2(\Omega)}+\va^{\f12}\|v_{xy}\|\nonumber\\
&&+\va^{\f12}\|\mathbf{g}\|_{H^1}\nonumber\\
\lesssim&&\va^{\f12}\|g_{0}\|_{H^2}+\|\mathbf{g}\|_{L^2}+|g_2(0,\cdot)|_{L^2} +\|g_0\|_{H^1}+\|\curl\mathbf{g}\|+\va^{\f12}\|\curl\mathbf{g}\|_{H^1}\nonumber\\
&&+\va|g_{0yy}(0,\cdot)|+\va^{\f12}\|\mathbf{g}\|_{H^1}+\va^{\f12}|g_{2y}(0,\cdot)|.\label{den.4}
\end{eqnarray}
combing (\ref{den.2})-(\ref{den.4}), we have
\begin{eqnarray}
&&\va^{\f12}\eta\|\rho_{xx}\|+\va^{\f12}\|\nabla\rho_{y}\|+\va\eta^2\sup_{x\in[0,L]}|\rho_{xx}(x,\cdot)|_{L^2}+\va\eta\sup_{x\in[0,L]}|\nabla\rho_{y}(x,\cdot)|L^2\nonumber\\
\lesssim&&\va^{\f12}\|g_{0}\|_{H^2}+\|\mathbf{g}\|_{L^2}+|g_2(0,\cdot)|_{L^2} +\|g_0\|_{H^1}+\|\curl\mathbf{g}\|+\va^{\f12}\|\curl\mathbf{g}\|_{H^1}\nonumber\\
&&+\va|g_{0yy}(0,\cdot)|+\va^{\f12}\|\mathbf{g}\|_{H^1}+\va^{\f12}|g_{2y}(0,\cdot)|,
\end{eqnarray}
and we have finished the proof of Lemma \ref{den.lem}.
\end{proof}
Combing the momentum equations (\ref{3.4})-(\ref{3.5}) and Lemma \ref{den.lem} we have the following Corollary:
\begin{corollary}\label{cor}
Let $(u,v,\rho)$ be solutions to system (\ref{3.3})-(\ref{boundary.2}), and the assumptions in Theorem \ref{th3.1} are satisfied. Moreover, we take $t=1$, then we have:
\begin{eqnarray*}
&&\va^{\f32}\|\nabla\dv\mau_y\|+\va^2\eta\|\dv\mau_{xx}\|\nonumber\\
\leq&&C(\va^{\f12}\|g_{0}\|_{H^2}+\|\mathbf{g}\|_{L^2}+|g_2(0,\cdot)|_{L^2} +\|g_0\|_{H^1}+\|\curl\mathbf{g}\|+\va^{\f12}\|\curl\mathbf{g}\|_{H^1}\nonumber\\
&&+\va|g_{0yy}(0,\cdot)|+\va^{\f12}\|\mathbf{g}\|_{H^1}+\va^{\f12}|g_{2y}(0,\cdot)|),
\end{eqnarray*}
here the constant $C$ is independent of $\va,\eta,\alpha.$

\end{corollary}
Finally, we have the following weighted estimates for $\nabla^3\mau$:
\begin{lemma}\label{u3.lem}
Let $(u,v,\rho)$ be solutions to the approximate system (\ref{app}), and the assumptions in Lemma \ref{lem.h1} are satisfied,  then we have:
\begin{eqnarray*}
&&\va^{\f32}\|(L-x)\nabla^2u_y\|+\va^2\eta\|(L-x)u_{xxx}\|+\va^{\f32}\|(L-x)\nabla^3v\|\nonumber\\
\leq&&C(\va^{\f12}\|g_{0}\|_{H^2}+\|\mathbf{g}\|_{L^2}+|g_2(0,\cdot)|_{L^2} +\|g_0\|_{H^1}+\|\curl\mathbf{g}\|+\va^{\f12}\|\curl\mathbf{g}\|_{H^1}\nonumber\\
&&+\va|g_{0yy}(0,\cdot)|+\va^{\f12}\|\mathbf{g}\|_{H^1}+\va^{\f12}|g_{2y}(0,\cdot)|),
\end{eqnarray*}
here the constant $C$ is independent of $\va,\eta,\alpha.$
\end{lemma}
\begin{proof}
First we  write the momentum equations in the following form:
\begin{equation}\label{u2}\begin{cases}
-\va\Delta u=-u_su_x-u_{sy}v+\va\partial_x\mathrm{div}\mathbf{u}+\gamma\rho_x
-g_1,\\
-\va\Delta v=-u_sv_x
+\va\partial_y\mathrm{div}\mathbf{u}-\gamma\rho_y+g_2.\end{cases}
\end{equation}
Then we differentiate $(\ref{u2})_1$ with respect to $y$, using boundary condition (\ref{boundary}) to have
\begin{equation}\label{u.1}\begin{cases}
-\va\Delta u_y=-(u_su_x-u_{sy}v)_y+\va\partial_{xy}\mathrm{div}\mathbf{u}+\gamma\rho_{xy}
-g_{1y},\\
u_y|_{y=0,2}=u_y|_{x=0}=0.\end{cases}
\end{equation}
Then we multiply (\ref{u.1}) with  $\va^2(L-x)^2u_{yyy}$ to have
\begin{eqnarray*}
&&\va^3\int_\Omega[(L-x)^2u_{yyy}^2+(L-x)^2u_{xyy}^2-2(L-x)u_{xyy}u_{yy}]dxdy\\
=&&\va^2\int_\Omega(L-x)^2u_{yyy}[-(u_su_x-u_{sy}v)_y+\va\partial_{xy}\mathrm{div}\mathbf{u}-\gamma\rho_{xy}+g_{1y}]dxdy,
\end{eqnarray*}
which combined with equation $(\ref{u.1})_1$, (\ref{app.2}), Lemma \ref{den.lem} and Corollary \ref{cor} implies
\begin{eqnarray}
&&\va^{\f32}\|(L-x)\nabla^2u_y\|\nonumber\\
\lesssim&& \va^{\f32}\|u_{yy}\|+\va^{\f12}[\|u_x\|_{H^1}+\|v\|_{H^1}+\va\|\partial_{xy}\dv\mau\|+\|\rho_{xy}\|+\|g_{1y}\|]\nonumber\\
\lesssim&&\va^{\f12}\|g_{0}\|_{H^2}+\|\mathbf{g}\|_{L^2}+|g_2(0,\cdot)|_{L^2} +\|g_0\|_{H^1}+\|\curl\mathbf{g}\|+\va^{\f12}\|\curl\mathbf{g}\|_{H^1}\nonumber\\
&&+\va|g_{0yy}(0,\cdot)|+\va^{\f12}\|\mathbf{g}\|_{H^1}+\va^{\f12}|g_{2y}(0,\cdot)|.\label{u.2}
\end{eqnarray}
Similarly  we differentiate $(\ref{u2})_2$ with respect to $x$, using boundary condition (\ref{boundary}) to have
\begin{equation}\label{v.1}\begin{cases}
-\va\Delta v_x=-u_sv_{xx}+\va\partial_{xy}\mathrm{div}\mathbf{u}+\gamma\rho_{xy}
-g_{2x},\\
v_x|_{y=0,2}=v_x|_{x=0}=0.\end{cases}
\end{equation}
Then we multiply $(\ref{v.1})_1$ with $\va^2(L-x)^2v_{xxx}$ to have
\begin{eqnarray*}
&&\va^3\int_\Omega[(L-x)^2v_{xxx}^2+(L-x)^2v_{xxy}^2-2(L-x)v_{xxy}v_{xy}]dxdy\\
=&&\va^2\int_\Omega(L-x)^2v_{xxx}[-u_sv_{xx}+\va\partial_{xy}\mathrm{div}\mathbf{u}-\gamma\rho_{xy}+g_{2x}]dxdy,
\end{eqnarray*}
which combined with equation $(\ref{u.1})_1$  (\ref{app.2}), Lemma \ref{den.lem} and Corollary \ref{cor} implies
\begin{eqnarray}
&&\va^{\f32}\|(L-x)\nabla^2v_x\|\leq \va^{\f12}[\|v\|_{H^1}+\va\|\partial_{xy}\dv\mau\|+\|\rho_{xy}\|+\|g_{2x}\|]\nonumber\\
\leq&&C(\va^{\f12}\|g_{0}\|_{H^2}+\|\mathbf{g}\|_{L^2}+|g_2(0,\cdot)|_{L^2} +\|g_0\|_{H^1}+\|\curl\mathbf{g}\|+\va^{\f12}\|\curl\mathbf{g}\|_{H^1}\nonumber\\
&&+\va|g_{0yy}(0,\cdot)|+\va^{\f12}\|\mathbf{g}\|_{H^1}+\va^{\f12}|g_{2y}(0,\cdot)|).\label{v.2}
\end{eqnarray}
Finally Lemma \ref{u3.lem} follows immediately by combing (\ref{u2}), (\ref{u.2}) and (\ref{v.2}).
\end{proof}
\subsection{Solutions to the linearized system}
In this section we will prove the existence and uniform-in-$\va$ estimates of solutions to system (\ref{linear.1})-(\ref{linear.3}) with boundary condition (\ref{boundary}).
The main result of this section reads as:
\begin{theorem}\label{main.linear}
Let $2\leq p<\infty$, for given $\mau^\va,\mathbf{g},g_0$ satisfying (\ref{in2}) and (\ref{in1}), there exists a unique triple $(u,v,\rho)\in$ satisfying system (\ref{linear.1})-(\ref{linear.3}) with boundary condition (\ref{boundary}) as well as the following estimates:
\begin{eqnarray}
&&\|\rho\|_{\mathcal{Y}}+\va\|u_y-v_x\|_{H^2}+\|\mau\|_{\mathcal{X}}\nonumber\\
\leq&&C(\|g_0\|_{H^1}+\va|g_{0y}(0,\cdot)|_{L^2}+\|\mathbf{g}\|_{L^2}+\|\curl\mathbf{g}\|_{L^2}+|g_2(0,\cdot)|_{L^2})
\label{th.3.1}
\end{eqnarray}
and
\begin{eqnarray}
&&\|\mathbf{u}\|_{\mathcal{B}}+\|\rho\|_{\mathcal{A}}\nonumber\\
\leq&& C(\va^{\f12}\|g_{0}\|_{H^2}+\|\mathbf{g}\|_{L^2}+|g_2(0,\cdot)|_{L^2} +\|g_0\|_{H^1}+\|\curl\mathbf{g}\|+\va^{\f12}\|\curl\mathbf{g}\|_{H^1}\nonumber\\
&&+\va|g_{0yy}(0,\cdot)|+\va^{\f12}\|\mathbf{g}\|_{H^1}+\va^{\f12}|g_{2y}(0,\cdot)|
\label{th.3}
\end{eqnarray}
where the norm $\|\cdot\|_{\mathcal{A}},\|\cdot\|_{\mathcal{B}},\|\cdot\|_{\mathcal{X}},\|\cdot\|_{\mathcal{Y}}$ are defined in (\ref{a}),(\ref{b}),(\ref{x}),(\ref{y}), and the constant $C$ is independent of $\va,\eta, L$.
\end{theorem}
\begin{proof}
As the estimates in Theorem \ref{thapp} and section 2.2 are independent of $\alpha$,  if we  take the limit $\alpha\rightarrow0$, then there exists a subsequence (still denoted by $(u^{app},v^{app},\rho^{app})$ and a unique triple $(u,v,\rho)\in W^{2,p}\times W^{2,p}\times H^2$ satisfying
\begin{eqnarray*}
u^{app}\rightarrow u\ \ \text{in}\ W^{2,p},\ \text{for}\ \forall\ 1<p<\infty;\\
v^{app}\rightarrow v\ \ \text{in}\ W^{2,p},\ \text{for}\ \forall\ 1<p<\infty;\\
\rho^{app}\rightarrow \rho\ \ \text{in}\ W^{1,p},\ \text{for}\ \forall\ 1<p<\infty,
\end{eqnarray*}
and estimates (\ref{th.3.1}) and (\ref{th.3}).
Besides, as
$$2\eta^2\rho_x^{app}(0,y)=\f1{ u^\va}\int_0^y[ (u_y^{app}-v_x^{app})_x^{\alpha}-(u_y^{app}-v_x^{app})_x](0,s)ds,$$
and we have the estimates:
\begin{eqnarray*}&&\va^{\f12}|\curl\mau_{x}^{app}(0,\cdot)|+\va^{\f32}|(\curl\mau^{app})_{xy}(0,\cdot)|\\
\leq&& C(\va^{\f12}\|g_{0}\|_{H^2}+\|\mathbf{g}\|_{L^2}+|g_2(0,\cdot)|_{L^2} +\|g_0\|_{H^1}+\|\curl\mathbf{g}\|+\va^{\f12}\|\curl\mathbf{g}\|_{H^1}\nonumber\\
&&+\va|g_{0yy}(0,\cdot)|+\va^{\f12}\|\mathbf{g}\|_{H^1}+\va^{\f12}|g_{2y}(0,\cdot)|),\end{eqnarray*}
here the constant $C$ is independent of $\alpha$. So if we  take the limit $\alpha\rightarrow0$, then we have
$$\rho_x|_{x=0}=0,$$
and we finish  the proof of Theorem \ref{main.linear}.
\end{proof}
%%%%%%%%%%%%%%%%%%%%%%%%%%%%
\renewcommand{\theequation}{\thesection.\arabic{equation}}
\setcounter{equation}{0}
\section{Solutions to the nonlinear system}

In this section we will prove the main result Theorem \ref{th.main}. As there is a term $(\mathbf{U}_0+\bar {\mathbf{u}}+\va^{\f12+}\mathbf{u}^{n})\cdot\nabla\rho^{n+1}$
on the left-hand side of the mass equation (\ref{3.1})$_1$, so that there will be a loss of regularity on the right-hand side of the equation.
It is not possible to directly show the convergence of the sequence $\{(\mathbf{u}^{n},\rho^{n})\}$ in $\mathcal{B}\times \mathcal{A}$ by a fixed point theory
  to obtain the strong solution of the nonlinear system (\ref{0.3})-(\ref{0.5}). Instead, we will establish the convergence by showing that $\{(\mathbf{u}^{n},\rho^{n})\}$ is a Cauchy
sequence in $\mathcal{X}\times \mathcal{Y}$ with a uniform upper bound in $\mathcal{B}\times \mathcal{A}$.
Then, by taking to the limit we obtain that the limit function $(\mathbf{u},\rho)$ is a solution to  system (\ref{0.3})-(\ref{0.5}) with estimate (\ref{th.3}).

The following Lemma gives the uniform bound of $\{(\mathbf{u}^{n},\rho^{n})\}$ in $\mathcal{B}\times \mathcal{A}$.

\begin{lemma}\label{lemma3.1}
Let $\{(\mathbf{u}^{n},\rho^{n})\}$ be a sequence of solutions to
the system (\ref{3.1}), then there exists a constant $M>0$, depending only on $p$, such that
for any $n\in Z^+$ we have
\begin{eqnarray}
&&\|\mathbf{u}^n\|_{\mathcal{B}}+\|\rho^n\|_{\mathcal{A}}\leq M\va^{\f\sigma2}\label{n1}.
\end{eqnarray}

\end{lemma}

\begin{proof}
First of all, by Theorem \ref{main.linear}  we have
\begin{eqnarray}
&&\|\mathbf{u}^{n+1}\|_{\mathcal{B}}+\|\rho^{n+1}\|_{\mathcal{A}}\nonumber\\
\leq&& C_1\{\va^{\f12}\|\curl\mathbf{g}(\mau^n,\rho^n)\|_{H^1}+\va^{\f12}\|g_0(\mau^n,\rho^n)\|_{H^2}+\|g_0(\mau^n,\rho^n)\|_{H^1}+\|\mathbf{g}(\mau^n,\rho^n)\|\nonumber\\
&&+\va|g_{0yy}(\mau^n,\rho^n)(0,\cdot)|+\va^{\f12}|g_{2y}(\mau^n,\rho^n)(0,\cdot)|+\va^{\f12}\|\mathbf{g}(\mau^n,\rho^n)\|_{H^1}\nonumber\\
&&+\|\curl\mathbf{g}(\mau^n,\rho^n)\|+|g_2(\mau^n,\rho^n)(0,\cdot)|\},\label{n2}
\end{eqnarray}
where the constant $C_1$ depends only on $\Omega$ and $ p$.

Now we prove (\ref{n1}) by an induction argument.

 When $n=1$, by the assumption that $(\mathbf{u}^0,\rho^0)=(0,0)$ we have
 $$g_0(0,0,0)=g_{0r},\ g_1(0,0,0)=g_{1r},\ g_2(0,0,0)=g_{2r}.$$
Recalling that $\eta<\va^{\f12+}$, direct computation shows
\begin{eqnarray*}
&&\|\mathbf{u}^{1}\|_{\mathcal{B}}+\|\rho^{1}\|_{\mathcal{A}}\nonumber\\
\leq&&C_1\{\va^{\f12}\|\curl\mathbf{g}(\mau^0,\rho^0)\|_{H^1}+\va^{\f12}\|g_0(\mau^0,\rho^0)\|_{H^2}+\|g_0(\mau^0,\rho^0)\|_{H^1}+\|\mathbf{g}(\mau^0,\rho^0)\|\nonumber\\
&&+\va|g_{0yy}(\mau^0,\rho^0)(0,\cdot)|+\va^{\f12}|g_{2y}(\mau^0,\rho^0)(0,\cdot)|+\va^{\f12}\|\mathbf{g}(\mau^0,\rho^0)\|_{H^1}\nonumber\\
&&+\|\curl\mathbf{g}(\mau^0,\rho^0)\|+|g_2(\mau^0,\rho^0)(0,\cdot)|\}\\
\leq&& C_1(\va^{\f12-}\alpha_2+\va^{-\f12-}\Lambda)\leq C_1\va^{\f{3\sigma}4}.
\end{eqnarray*}

Taking $M\geq \va^{\f\sigma4}C_1$, then (\ref{n1}) holds for $n=1$.

 Now, assuming that for any $1\leq k\leq n$,
\begin{eqnarray}
\|\mathbf{u}^k\|_{\mathcal{B}}+\|\rho^k\|_{\mathcal{A}}\leq M\va^{\f\sigma2}\label{n4}
\end{eqnarray}
and denoting
$$A_k=\|\mathbf{u}^k\|_{\mathcal{B}}+\|\rho^k\|_{\mathcal{A}},$$
First of all, by Sobolev imbedding theorem, we have for any $n\geq1$
\begin{eqnarray}
&&\va^+|\rho^n|_{L^\infty}+\va^{\f12+}|\rho^n_y|_{L^\infty}+\va\eta|\rho_x^n|_{L^\infty}+\va^+|\mau^n|_{L^\infty}+\va^{\f12+}|\nabla\mau^n|_{L^\infty}
+\va|v^n_{yy}(0,\cdot)|\nonumber\\
\lesssim&&\|\rho^n\|_{H^1}+\va\eta\|\rho^n_{xx}\|+\va^{\f12}\|\nabla\rho^n_{y}\|
+\va\eta|\nabla\rho^n_{y}|_{L^\infty_xL^2_y}+\|\mau^n\|_{H^1}+\va^{\f12}\|\nabla^2\mau^n\|\nonumber\\
&&+\va^{2-\f2p}\|\mau^n\|_{2,p}+\va^{\f32}\|(L-x)\nabla^3v^n\|)\nonumber\\
\lesssim&&\|\mathbf{u}^n\|_{\mathcal{B}}+\|\rho^n\|_{\mathcal{A}}=A_n.\label{3.6}
\end{eqnarray}
Then by (\ref{3.6}) and direct computation
\begin{eqnarray*}
&&\|g_0(\mau^n,\rho^n)\|_{H^1}+\va^{\f12}\|g_0(\mau^n,\rho^n)\|_{H^2}+|g_{0yy}(\mau^n,\rho^n)(0,\cdot)|\\
\lesssim&& \va\eta^2\{\|\rho^n\dv\mau^n\|_{H^1}+\va^{\f12}\|\rho^n\dv\mau^n\|_{H^2}+\va^{\f12}\|\dv\mau^n\|_{H^2}+\|\mau^n\|_{H^2}+\|\rho^n\|_{H^2}\\
&&+|(\rho^n\dv\mau^n)_{yy}(0,\cdot)|+L|\dv\mau^n(0,\cdot)|_{H^2}+|\rho^n(0,\cdot)|_{H^2}\}+A_1\\
\lesssim&&A_n^2+\eta A_n+A_1.
\end{eqnarray*}
Next, Since
\begin{eqnarray}
&&\partial_yg_{11}(\rho,u,v)-\partial_xg_{21}(\rho,u,v)\nonumber\\
=&&\partial_y\left(\{(2\alpha_2 -\gamma\rho_x)[(\rho^\va)^{\gamma-1}-1]-\va^{-\f12-}\gamma\bar\rho_x(\rho^\va)^{\gamma-1}\}
-\eta^2\rho(u^\va u^\va_x+v^\va u^\va_y)\right)\nonumber\\
&&-\partial_x\left(\{\gamma\rho_y[(\rho^\va)^{\gamma-1}-1]-\va^{-\f12-}\gamma\bar\rho_y(\rho^\va)^{\gamma-1}\}-\eta^2 \rho(u^\va v^\va_x+v^\va v^\va_y)\right)\nonumber\\
=&&(2\alpha_2 -\gamma\rho_x)(\rho^\va)_y^{\gamma-1}-\gamma\rho_y(\rho^\va)_x^{\gamma-1}-\partial_y[\va^{-\f12-}\gamma\bar\rho_x(\rho^\va)^{\gamma-1}]\nonumber\\
&&-\partial_y[\eta^2\rho(u^\va u^\va_x+v^\va u^\va_y)]+\partial_x[\va^{-\f12-}\gamma\bar\rho_y(\rho^\va)^{\gamma-1}+\eta^2 \rho(u^\va v^\va_x+v^\va v^\va_y)],\label{curlg}
\end{eqnarray}
while
$$(\rho^\va)_x^{\gamma-1}\sim\eta^2(-\va\f2\gamma\alpha_2+\bar\rho_x+\va^{\f12+}\rho_x),\ (\rho^\va)_y^{\gamma-1}\sim\eta^2(\bar\rho_y+\va^{\f12+}\rho_y),$$
and
\begin{eqnarray*}
&&\partial_yg_{12}(u,v)-\partial_xg_{22}(u,v)\nonumber\\
=&&-\va^{\f12+} u(u_y-v_x)_x-\va^{\f12+} v(u_y-v_x)_y-\bar u(u_y-v_x)_x-\bar v(u_y-v_x)_y\nonumber\\
&&-\eta^2(-\f2\gamma\alpha_2\va x+\bar\rho)[u^\va(u_y-v_x)_x+v(u^\va_y-v^\va_x)_y+\bar v(u_y-v_x)_y]+\va^{\f12+} h(\nabla\mau,\mau),
\end{eqnarray*}
here $h(\nabla\mau,\mau)$ is a known function which could be well controlled.
So if we denote by $\mau^{\va,n}=(u_s+ \bar u+\va^{\f12+} u^{n}),\bar v+\va^{\f12+} v^{n})$,  we have by (\ref{3.6})
\begin{eqnarray*}
&&\va^{\f12}\|\curl\mathbf{g}(\mau^n,\rho^n)\|_{H^1}+\va^{\f12}\|\mathbf{g}(\mau^n,\rho^n)\|_{H^1}+\|\mathbf{g}(\mau^0,\rho^0)\|\nonumber\\
&&+\va^{\f12}|g_{2y}(\mau^0,\rho^0)(0,\cdot)|+\|\curl\mathbf{g}(\mau^0,\rho^0)\|+|g_2(\mau^0,\rho^0)(0,\cdot)|\\
\lesssim&&\va^{1+}\eta^2\|\rho^n_x\rho^n_y\|_{H^1}+\va\eta^2\|\nabla\rho\|_{H^1}+\eta^2\va^{1+}|\rho^n|_{L^\infty}\|\mau^{n}\|_{H^1}^2\\
&&+\va^{1+}(1+\eta^2|\rho^n|_{L^\infty})|\mau^n|_{L^\infty}\|u_y^n-v^n_x\|_{H^2}+\va^{1+}\eta^2|\rho^n|_{L^\infty}\|\rho^n\|_{H^2}\\
&&+\va^{1+}|\mau^n|_{L^\infty}\|\mau^n\|_{H^2}+\va\eta^2\|\rho^n\|_{H^1}+\va^{1+}\|\mau^n\|_{H^2}+\va^{1+}\eta^2|\rho^n|_{L^\infty}|\rho_{yy}(0,\cdot)|\\
&&+\va\eta^2|\rho_y(0,\cdot)|[|\rho_y(0,\cdot)|+|v^{\va,n}v^{\va,n}_y|_{L^\infty}]+\va^{1+}(1+|v^n|_{L^\infty})|v^n(0,\cdot)|_{H^2}\\
\lesssim&&A_n^2+\va^{\f12} A_n
\end{eqnarray*}
Combing above, we have
\begin{eqnarray*}
A_{n+1}\leq C(A_n^2+\va^{\f12}A_n+A_1)
\end{eqnarray*}
Therefore, (\ref{n1}) follows immediately.
\end{proof}

Based on Lemma \ref{lemma3.1}, we can prove that $\{(\mathbf{u}^{n},\rho^{n})\}$ is a Cauchy sequence in $\mathcal{X}\times \mathcal{Y}$. More precisely, we have the following Lemma:
\begin{lemma}\label{cs}
Let $\{(\mathbf{u}^{n},\rho^{n})\}_{n=1}^\infty$ be the sequence of solutions to system (\ref{3.1}) with initial data $(\mathbf{u}^{0},\rho^{0})=(0,0,0)$, then we have
\begin{eqnarray}
\|\mathbf{u}^{n+1}-\mathbf{u}^{n}\|_{\mathcal{X}}+\|\rho^{n+1}-\rho^{n}\|_{\mathcal{Y}}\leq
\f12 \Big( \|\mathbf{u}^{n}-\mathbf{u}^{n-1}\|_{\mathcal{X}}+\|\rho^{n}-\rho^{n-1}\|_{\mathcal{Y}}\Big).
\end{eqnarray}
\end{lemma}

\begin{proof} A straightforward calculation gives
\begin{eqnarray}
&& \dv(\mau^{n+1}-\mau^n)+\eta^2\mau^{\va,n}\cdot\nabla(\rho^{n+1}-\rho^n)\nonumber\\
=&&g_{01}(\mathbf{u}^{n},\rho^{n})-g_{01}(\mathbf{u}^{n-1},\rho^{n-1})
- \eta^2\va^{\f12+}(\mathbf{u}^{n}-\mathbf{u}^{n-1})\cdot\nabla\rho^{n}\nonumber\\
\triangleq&&\hat g_0\qquad\text{in}~\Omega,\nonumber\\
&&\  u_s(u^{n+1}-u^n)_x+u_{sy}(v^{n+1}-v^n)-\va\Delta (u^{n+1}-u^n)
-\va\partial_x\mathrm{div}(\mau^{n+1}-\mau^n)\nonumber\\
&&\quad+\gamma(\rho^{n+1}-\rho^n)_x\nonumber\\
=&&g_{11}(\mathbf{u}^{n},\rho^{n})+g_{12}(\mathbf{u}^{n},\rho^{n})-g_{11}(\mathbf{u}^{n-1},\rho^{n-1})-g_{12}(\mathbf{u}^{n-1},\rho^{n-1})\nonumber\\
\triangleq&& \hat g_1,\qquad\text{in}~\Omega,\nonumber\\
&&\  u_s(v^{n+1}-v^n)_x-\va\Delta (v^{n+1}-v^n)
-\va\partial_y\mathrm{div}(\mau^{n+1}-\mau^n)+\gamma(\rho^{n+1}-\rho^n)_y\nonumber\\
=&&g_{21}(\mathbf{u}^{n},\rho^{n})+g_{22}(\mathbf{u}^{n},\rho^{n})-g_{21}(\mathbf{u}^{n-1},\rho^{n-1})-g_{22}(\mathbf{u}^{n-1},\rho^{n-1})\nonumber\\
\triangleq&&\hat g_2\qquad\text{in}~\Omega,
\end{eqnarray}
with boundary condition
\begin{equation}\begin{cases}
(\rho^{n+1}-\rho^{n})_x|_{x=0}=0,\\
(u^{n+1}-u^n)|_{x=0}=(v^{n+1}-v^n)_x|_{x=0}=0,\nonumber\\
(v^{n+1}-v^n)|_{x=L}=\curl(\mau^{n+1}-\mau^n)_x|_{x=L}=0,\\
(u^{n+1}-u^n)_y|_{y=0,2}=(v^{n+1}-v^n)|_{y=0,2}=0
\end{cases}\end{equation}
Next, by (\ref{3.6}) and direct computation we have
\begin{eqnarray*}
&&\|\hat g_0\|_{H^1}+|\hat g_{0y}(0,\cdot)|_{L^2}\\
\lesssim&&\va^{\f12+}\eta^2\|\rho^n\dv\mau^n-\rho^{n-1}\dv\mau^{n-1}\|_{H^1}+\va^{\f12+}\eta^2\|(\mathbf{u}^{n}-\mathbf{u}^{n-1})\cdot\nabla\rho^{n}\|_{H^1}\\
&&+\va^{\f12+}\eta^2|(\rho^n\dv\mau^n-\rho^{n-1}\dv\mau^{n-1})(0,\cdot)|_{H^1}+\va^{\f12+}\eta^2|(\mathbf{u}^{n}-\mathbf{u}^{n-1})\cdot\nabla\rho^{n}(0,\cdot)|_{H^1}\\
&&+\eta^2\va^{\f12+}(\|\mau^n-\mau^{n-1}\|_{H^2}+\|\dv\mau^n-\dv\mau^{n-1}\|_{H^1}+\|\rho^n-\rho^{n-1}\|_{H^1})\\
&&+\eta^2\va^{\f12+}(|(\dv\mau^n-\dv\mau^{n-1})(0,\cdot)|_{H^1}+|(\rho^n_y-\rho_y^{n-1})(0,\cdot)|_{L^2})\\
\lesssim&&  \va^{\f12}(\|\mathbf{u}^{n+1}-\mathbf{u}^{n}\|_{\mathcal{X}}+\|\rho^n-\rho^{n-1}\|_{\mathcal{Y}})
\end{eqnarray*}
and using (\ref{3.6}) and (\ref{curlg}) we obtain
\begin{eqnarray*}
&&\|\hat {\mathbf{g}}\|+|\hat g_2(0,\cdot)|+\|\curl\hat{\mathbf{g}}\|\\
\lesssim&&\va^{\f12} \|\mathbf{u}^{n}-\mathbf{u}^{n-1}\|_{\mathcal{X}}+\eta\|\rho^{n}-\rho^{n-1}\|_{\mathcal{Y}}+\|(\bar vv_y^n-\bar vv_y^{n-1})(0,\cdot)\|\\
&&+\va^{\f12+}\|(v^nv_y^n-v^{n-1}v_y^{n-1})(0,\cdot)\|+\va^{\f12+}\|\mau^n\nabla\curl\mau^n-\mau^{n-1}\nabla\curl\mau^{n-1}\|\\
\lesssim&&\va^{+}\|\mau^n-\mau^{n-1}\|_{\mathcal{X}}+\eta\|\rho^n-\rho^{n-1}\|_{\mathcal{Y}}
\end{eqnarray*}
Combing above, using  (\ref{th.3.1}) in Theorem \ref{main.linear} and the smallness of $\va$ we obtain
\begin{eqnarray}
&&\|\mau^n-\mau^{n-1}\|_{\mathcal{X}}+\|\rho^n-\rho^{n-1}\|_{\mathcal{Y}}\nonumber\\
\lesssim &&\|\hat g_0\|_{H^1}+|\hat g_{0y}(0,\cdot)|_{L^2}+\|\hat {\mathbf{g}}\|+|\hat g_2(0,\cdot)|+\|\curl\hat{\mathbf{g}}\|\nonumber\\
\leq&& \f12\|\mau^n-\mau^{n-1}\|_{\mathcal{X}}+\f12\|\rho^n-\rho^{n-1}\|_{\mathcal{Y}}
\end{eqnarray}
and we complete the proof.
\end{proof}

\textbf{Proof of Theorem \ref{th.main}.} \ Lemma \ref{cs} implies that
$\{(\mathbf{u}^{n},\rho^{n})\}$ is a Cauchy sequence in
$\mathcal{X}\times \mathcal{Y}$. Hence, with the help of Lemma \ref{lemma3.1} we see that
there exists a unique $({\mathbf{u}},\rho)\in \mathcal{B}\times \mathcal{A}$,
 such that for any $1<p'<p$,
$$\mathbf{u}^{n}\rightarrow\mathbf{u}\;\quad\text{strongly in }\; W^{2,p'}(\Omega), $$
$$\rho^{n}\rightarrow\rho\;\quad\text{strongly in }\; W^{1,p'}(\Omega).$$

 Taking limit on both sides of (\ref{3.1}), we conclude that
 $(\mathbf{u},\rho)$ satisfies the nonlinear system (\ref{0.3})-(\ref{0.5}) with boundary condition (\ref{0.b1})-(\ref{0.b2}). Moreover, the estimate (\ref{main.1}) follows from (\ref{n1})  immediately and (\ref{main.2})-(\ref{main.3}) follows immediately from Sobolev imbedding Theorem.

Next, if $(\mathbf{u},\rho)$ and $(\hat{\mathbf{u}},\hat{\rho})$ are
two solutions of the system (\ref{1.1}), then by a process similar to that used
in Lemma \ref{3.5}, we infer that
\begin{eqnarray}
\|\mathbf{u}-\hat{\mathbf{u}}\|_{1,2;\Omega}+\|\rho-\hat{\rho}\|_{2;\Omega}\leq
\f12(\|\mathbf{u}-\hat{\mathbf{u}}\|_{1,2;\Omega}+\|\rho-\hat{\rho}\|_{2;\Omega}).
\end{eqnarray}
Consequently, we obtain the uniqueness of the solutions. This completes the proof of Theorem 1.1.

%%%%%%%%%%%%%%%%%%%%%%%%%%%%


\begin{thebibliography}{50}

\bibitem{ADM}Antonelli, P.; Dolce, M.; Marcati, P.: Linear stability analysis of the homogeneous Couette flow in a 2D isentropic compressible fluid
Ann. PDE 7 (2021), no. 2, Paper No. 24, 53 pp.

\bibitem{B-N} J. B\v{r}ezina and A.Novotn\'y: On weak solutions of steady Navier-Stokes equations for
monoatomic gas, {\it Comment. Math. Univ. Carolinae} {\bf 49} no.
4,(2008), 611-632.

\bibitem{C-J} H. Choe and B. Jin: Existence of solutions of stationary compressible Navier-Stokes Equations
with large force, {\it J. Funct. Anal.} {\bf 177} (2000), 54-88.

\bibitem{D-J-J-Y} C. Dou, F. Jiang, S. Jiang, Y. Yang: Existence of strong solutions to the steady
Navier-Stokes equations for a compressible heat-conductive fluid with large forces,
{\it J. Math. Pures Appl.} {\bf 103} (2015), 1163-1197.

\bibitem{G-T} D. Gilbarg, N.S. Trudinger: {\it Elliptic Partial Differential Equations of Second Order},
2nd Edition, Springer-Verlag, 1983.
\bibitem{GJZ}
Guo, Yan; Jiang, Song; Zhou, Chunhui: Steady viscous compressible channel flows,
{\it SIAM J. Math. Anal.} 47  no. 5, 3648-3670, 2015

\bibitem{LW}Liu, C.J., Wang, Y.G.: Stability of boundary layers for the nonisentropic compressible
circularly symmetric 2D flow. {\it SIAM J. Math. Anal}. 46(1), 256-309, 2014



\bibitem{K}Kagei, Y. Asymptotic behavior of solutions to the compressible Navier-Stokes equation around a parallel flow
Arch. Ration. Mech. Anal. 205 (2012), no. 2, 585-650.

\bibitem{J-K}C.H. Jun, J.R. Kweon, For the stationary compressible viscous Navier-Stokes equations with no-slip condition
on a convex polygon. {\it J. Diff. Eqns.} {\bf 250} (2011), 2440-2461.

\bibitem{K-K}R.B. Kellogg, J.R. Kweon,compressible Navier-Stokes
equations in a bounded domain with inflow boundary condition, {\it SIAM J. Math. Anal.} \textbf{28} (1997), 94-108.

%\bibitem{k-k-1}J. R. Kweon,R. B. Kellogg, Smooth solution of the compressible Navier-Stokes
%Equations in an unbounded domain with inflow boundary condition, J. Math. Anal.
%App. 220 (1998), 657¨C675.

\bibitem{JK} Jerison, D.; Kenig, C.: Inhomogeneous Dirichlet problem in Lipschitz domains. J. Funct. Anal. 130 (1995), no. 1, 161-219.

\bibitem{J-Zhou} Song Jiang, Chunhui Zhou: Existence of weak
solutions to the three-dimensional steady compressible Navier-Stokes
equations, {\it Ann.I.H. Poincar\'{e} Analyse non lin\'{e}aire} {\bf
28} (2011), 485-498

\bibitem{k-k-2}J.R. Kweon, R.B. Kellogg, Regularity of solutions to the Navier-Stokes equations for compressible barotropic
flows on a polygon, {\it Arch. Rational Mech. Anal.} {\bf 163} (2002), 35-64.

%\bibitem{k}J.R. Kweon, A jump discontinuity of compressible viscous flows grazing a non-convex corner,
%{\it J. Math. Pures Appl.} {\bf 100} (2013), 410-432.

 \bibitem{K-K-3}J.R. Kweon, R.B. Kellogg, Regularity of solutions to the Navier-Stokes system for compressible flows
 on a polygon, {\it SIAM J. Math. Anal.} {\bf 35} (2004), 1451-1485.

 \bibitem{L-Z}Li, H.L.; Zhang, X.W.  Stability of plane Couette flow for the compressible Navier-Stokes equations with Navier-slip boundary
J. Differential Equations 263 (2017), no. 2, 1160¨C1187.

\bibitem{L-2} P.L. Lions: {\it Mathematical Topics in Fluid Mechanics, Vol.II, Compressible Models}.
Clarendon Press, Oxford, 1998.

 \bibitem{MWWZ}Masmoudi, N.; Wang, Y.X.; Wu, D.; Zhang, Z.F. Tollmien-Schlichting waves in the subsonic regime
Proc. Lond. Math. Soc. (3) 128 (2024), no. 3, Paper No. e12588, 112 pp

\bibitem{N-N-P} S.A. Nazarov, A. Novotn\'y, K. Pileckas, On steady compressible Navier-Stokes equations in plane domains with corners,
{\it Math. Ann}. {\bf 304} (1996), 121-150.

\bibitem{P} Paddick, M.: The strong inviscid limit of the isentropic compressible Navier-Stokes
equations with Navier boundary conditions. Discrete Contin. Dyn. Syst. 36(5), 2673-
2709, 2016

\bibitem{Pa}M. Padula, Existence and uniqueness for viscous steady
compressible motions, {\it Arch. Rational Mech. Anal}. \textbf{77} (1987), 89-102.

\bibitem{PT} T. Piasecki, On an inhomogeneous slip-inflow boundary
value problem for a steady flow of a viscous compressible fluid in a cylindrical domain,
{\it J. Diff. Eqns.} \textbf{248} (2010), 2171-2198.

\bibitem{P-P}T. Piasecki, M. Pokorn\'y, Strong solutions to the
Navier-Stokes-Fourier system with slip-inflow boundary conditions,
{\it Z. Angew. Math. Mech.} {\bf 94} (2014), 1035-1057.

\bibitem{P-R-S-1} P.I. Plotnikov, E.V. Ruban, J. Sokolowski, Inhomogeneous boundary
value problems for compressible Navier-Stokes an transport
equations, {\it J. Math. Pures Appl.} \textbf{92} (2009), 113-162.

\bibitem{P-R-S}P.I. Plotnikov, E.V. Ruban, J. Sokolowski, Inhomogeneous
boundary value problems for compressible Navier-Stokes Equations:
well-posedness and sensitivity analysis, {\it SIAM J. Math. Anal.} \textbf{40} (2008), 1152-1200.

\bibitem{P-W}Plotnikov, P. I.; Weigant, W.P.I. Plotnikov, Steady 3D viscous compressible flows with adiabatic exponent  $\gamma\in\ (1,\infty)$
J. Math. Pures Appl. (9) 104 (2015), no. 1, 58-82.


\bibitem{V} A. Valli, On the existence of stationary solutions to
compressible Navier-Stokes equations, {\it Ann. Inst. H. Poincar\'{e}} \textbf{4} (1987), 99-113.

\bibitem{V-Z}A. Valli, W.M. Zajaczkowski, Navier-Stokes equations for
compressible fluids: global existence and qualitabive properties of the solutions in the general case,
{\it Comm. Math. Phys.} \textbf{103} (1986), 259-296.

\bibitem{YZ1}Yang, A.; Zhang,Z.: Long-time instability of planar Poiseuille-type flow in compressible fluid. arxiv:2309.06775v1.

\bibitem{YZ}Yang, T.; Zhang, Z.: Linear instability analysis on compressible Navier-Stokes equations with strong boundary layer.
Arch. Ration. Mech. Anal. 247  no. 5, Paper No. 83, 53 pp 2023,



\bibitem{W}Wang, Y.: Uniform regularity and vanishing viscosity limit for the full compressible
Navier-Stokes system in three dimensional bounded domain. Arch. Ration. Mech. Anal.
221, 1345-1415, 2016


\bibitem{WWZ} Wang,C.,Wang,Y.X.,Zhang,Z.F.: Zero-viscosity limit of the compressible Navier-Stokes equations in the analytic setting. arxiv:2305.09393v1


\bibitem{WW}Wang, Y.G., Williams, M.: The inviscid limit and stability of characteristic boundary layers for the compressible Navier-Stokes equations with Navier-friction boundary conditions. Ann. Inst. Fourier (Grenoble) 62(6), 2257-2314, 2013

\bibitem{WXY}    Wang, Y., Xin, Z., Yong, Y.: Uniform regularity and vanishing viscosity limit for
the compressible Navier-Stokes with general Navier-slip boundary conditions in three
dimensional domains. SIAM J. Math. Anal. 47(6), 4123-4191, 2015

\bibitem{XY}Xin, Z., Yanagisawa, T.: Zero-viscosity limit of the linearized Navier-Stokes equations
for a compressible viscous fluid in the half plane. Commun. Pure Appl. Math. 52(4),
479-541, 1999

\bibitem{ZZ}Zeng, L., Zhang, Z., Zi, R.: Linear stability of the Couette flow in three dimensional
isentropic compressible Navier-Stokes equations. SIAM J. Math. Anal. 54(5), 5698-
5741, 2022

\end{thebibliography}
\end{document}